\documentclass[aop,MSNbibl,seceqn,citesort,dvips]{arximspdf}
\usepackage{graphicx}

\doi{10.1214/10-AOP640}
\volume{40}
\issue{3}
\pubyear{2012}
\firstpage{1009}
\lastpage{1040}

\makeatletter

\newcommand{\mod}{\operatorname{mod}}

\newcommand{\real}{{\mathbb R}}

\newcommand{\Dom}{\operatorname{Dom}}

\newcommand{\inte}{\mathbb{N}}

\newcommand{\trace}{\operatorname{trace}}

\newcommand{\dee}{\mathbb{D}}
\newcommand{\lee}{\mathbb{L}}

\newtheorem{prop}{Proposition}[section]
\newtheorem{theorem}[prop]{Theorem}
\newtheorem{lemma}[prop]{Lemma}
\newtheorem{corollary}[prop]{Corollary}

\makeatother

\begin{document}
\begin{frontmatter}

\title{Girsanov identities for Poisson measures under quasi-nilpotent transformations}
\runtitle{Girsanov identities for Poisson measures}

\begin{aug}
\author[A]{\fnms{Nicolas} \snm{Privault}\corref{}\thanksref{t1}\ead[label=e1]{nprivault@ntu.edu.sg}}
\runauthor{N. Privault}
\affiliation{Nanyang Technological University}
\address[A]{Division of Mathematical Sciences\\
School of Physical\\
\quad and Mathematical Sciences\\
Nanyang Technological University\\
SPMS-MAS-05-43, 21 Nanyang Link\\
Singapore 637371\\
\printead{e1}} 
\end{aug}

\thankstext{t1}{Supported by the Grant GRF 102309 from the Research
Grants Council of the Hong Kong Special Administrative Region of the
People's Republic of China.}

\received{\smonth{5} \syear{2010}}
\revised{\smonth{12} \syear{2010}}

%
\begin{abstract}
We prove a Girsanov identity on the Poisson space
for anticipating transformations that satisfy
a strong quasi-nilpotence condition.
Applications are given to the Girsanov theorem
and to the invariance of Poisson measures under random transformations.
The proofs use combinatorial identities for the central moments
of Poisson stochastic integrals.
\end{abstract}

%

\end{frontmatter}

\section{Introduction}\label{sec1}
The Wiener and Poisson measures are well known to be
quasi-invariant under adapted shifts.
This quasi-invariance property has been extended to anticipative
shifts by several authors; cf. \cite{ramer,kusuoka} and
\cite{uzbk} and references
therein in the Wiener case, and, for example,
\cite{smorodina1,girpri,girunif,priqi},
in the Poisson case.

In the anticipative case the corresponding
Radon--Nikodym density is usually written
as the product
\[
| {\det}_2 (I+\nabla u) | \exp \bigl( -\delta(u) - \tfrac{1}{2} \Vert u
\Vert^2 \bigr)
\]
of a Skorohod--Dol\'{e}ans exponential
with the Carleman--Fredholm determinant of the Malliavin
gradient $\nabla u$ of the shift $u$; cf.
\cite{ramer,kusuoka,uzbk}.
A similar formula can be obtained for
Poisson random measures; cf. Section \ref{s6}.

It has been noted in \cite{zeitouni} that the standard Dol\'{e}ans
form of the density for anticipative shifts $u\dvtx W\to H$ on the
Wiener space $W$ with Cameron--Martin space $H$
can be conserved [i.e., the Carleman--Fredholm determinant
${\det}_2 (I+\nabla u)$ equals one] when the gradient
$\nabla u$ of the shift $u$ is quasi-nilpotent,
that is,
%
%
\begin{equation}
\label{ee}
\lim_{n\to\infty} \Vert( \nabla u )^n \Vert^{1 / n}_{HS} = 0
\quad\mbox{or equivalently}\quad
\trace( \nabla u )^n = 0,\qquad
n\geq2;\hspace*{-35pt}
\end{equation}
cf. \cite{zeitouni} or Theorem 3.6.1 of \cite{uzbk}.
In particular, when $\nabla u$ is quasi-nilpotent and $\Vert u\Vert$
is constant, it has been shown in \cite{ustrotations} that
$\delta(u)$ has a centered Gaussian law with variance
$\Vert u\Vert^2$; cf. \cite{primoment} for a simplified proof.

In this paper we consider
the Poisson space $\Omega^X$ over a metric space $X$ with
$\sigma$-finite intensity measure $\sigma(dx)$, and
investigate the quasi-invariance of random
transformations $\tau( \omega, \cdot)$
which are assumed to be
quasi-nilpotent in the sense that the finite difference gradient
$D_s \tau( \omega, t )$ satisfies the cyclic finite difference
condition (\ref{cyclic3}) below, which is a strenghtened version of
(\ref{ee}).
We show in particular that such anticipating quasi-nilpotent
transformations are quasi-invariant, and their Radon--Nikodym
densities are given by Dol\'{e}ans stochastic exponentials with jumps.
This also extends and recovers other results on the invariance
of random transformations of Poisson measures; cf.~\cite{prinv}.

Our starting point is the classical Girsanov identity for Poisson
random measures which states that
%
%
\begin{equation}
\label{1}
E_\sigma\biggl[ \exp \biggl( - \int_X g(x) \sigma(dx) \biggr)
\prod_{x\in\omega} \bigl( 1 + g (x) \bigr) \biggr] = 1 ,
\end{equation}
and rewrites when $g = {\mathbf1}_A$ as
\[
E_\sigma\bigl[ e^{- r \sigma(A) } ( 1 + r )^{\omega(A)} \bigr] = 1,\qquad r \in\real,
\]
which is equivalent to the vanishing of the
expectation
\[
E [ C_n ( Z , \lambda) ] = 0,\qquad n \geq1 ,
\]
for $Z = \omega(A)$ a Poisson random variable with
intensity $\lambda= \sigma( A )$, where $C_n ( x , \lambda)$
is the Charlier polynomials of degree $n\in\inte$,
with generating function
\[
e^{- r \lambda} ( 1 + r )^x = \sum_{n=0}^\infty \frac{r^n}{n!} C_n ( x
, \lambda) ,\qquad r > -1 .
\]
It is well known, however, that $Z$ need not have
a Poisson distribution for $E [ C_n ( Z , \lambda) ]$
to vanish when $\lambda$ is allowed to be random.
Indeed, such an identity also holds in the random adapted case
under the form
%
%
\begin{equation}
\label{1ad}
E \bigl[ C_n \bigl( N_{\tau^{-1} ( t)} , \tau^{-1} ( t) \bigr) \bigr] = 0,
\qquad n \geq1,
\end{equation}
where $(N_t)_{t\in\real_+}$ is a standard Poisson process
generating a filtration $(\mathcal{F}_t)_{t\in\real_+}$
and $\tau(t)$ is an $\mathcal{F}_t$-adapted time change,
due to the fact that
\begin{eqnarray*}
&&C_n \bigl( N_{\tau^{-1} ( t )} , \tau^{-1} ( t ) \bigr)\\
&&\qquad = n! \int_0^\infty
\int_0^{t_n} \cdots \int_0^{t_2} d\bigl(N_{\tau^{-1} ( t_1) }- d\tau^{-1}
(t_1 ) \bigr) \cdots d\bigl( N_{\tau^{-1} ( t_n ) } - d \tau^{-1} ( t_n )
\bigr)\nonumber
\end{eqnarray*}
is an adapted $n$th order iterated multiple stochastic integral
with respect to the compensated point martingale
$(N_{\tau^{-1} ( t )} - \tau^{-1} ( t ) )_{t\in\real_+}$;\vadjust{\goodbreak} cf. \cite
{kailath} and~\cite{coursm}, page 320.
In this case we also have
\[
E_\sigma\bigl[ e^{ - r \tau^{-1} (t) } ( 1 + r )^{N_{\tau^{-1} ( t )}
}\bigr]=1,\qquad
r \in\real,
\]
and more generally
%
%
\begin{equation}
\label{asdf}
E_\sigma\biggl[ \exp \biggl( - \int_0^\infty g ( \tau( s ) )\, d s \biggr)
\mathop{\prod_{\Delta N_s =1}}_{0 < s < \infty} \bigl( 1 + g ( \tau( s ) ) \bigr) \biggr] = 1,
\end{equation}
under a Novikov-type integrability condition on $g\dvtx\real\to\real$;
cf., for example,~\cite{lepinglememin}.

In Corollary \ref{cjds} below we will extend the Girsanov
identity (\ref{asdf}) to random anticipating
processes indexed by an abstract space $X$,
by computing the expectation
\[
E_\sigma[ C_n ( \omega( A ) , \sigma( A ) )],\qquad
n\geq1,
\]
of the random Charlier polynomial
$C_n ( \omega( A ) , \sigma( A ) )$,
where $A (\omega)$ is a random, possibly anticipating set.
In particular we provide conditions on $A (\omega)$
for the expectation
$E_\sigma[C_n ( \omega( A ) , \sigma( A ) )]$, $n\geq1$,
to vanish; cf. Proposition \ref{p12} below.
Such conditions are satisfied, in particular, under the
quasi-nilpotence condition (\ref{cyclic3}) below
and include the adaptedness
of $(\tau( t ))_{t\in\real_+}$ above,
which recovers the classical adapted Girsanov identity
(\ref{asdf}) as a particular case; cf. Proposition \ref{cc11}.
As a consequence we will obtain a Girsanov theorem for random
transformations of Poisson samples on an arbitrary measure
space.

The above results will be proved using the Skorohod integral
and integration by parts on the Poisson space.
This type of argument has been applied in~\cite{prinv}
to the inductive computation of moments of Poisson stochastic
integrals and to the invariance of the Skorohod integral under
random intensity preserving transformations.
However, the case of Charlier polynomials
is more complicated, and it leads to Girsanov identities
and a Girsanov theorem as additional applications.

Since our use of integration by parts formulas
and moment identities relies on compensated Poisson
stochastic integrals, we will need to work with a~family
$B_n ( y , \lambda)$ of polynomials such that
\[
B_n ( y , - \lambda) = E_\lambda[ ( Z + y - \lambda)^n ] ,
\]
where $Z$ is a Poisson random variable with intensity $\lambda>0$,
and which are related to the Charlier polynomials by the relation
\[
C_n ( y , \lambda) = \sum_{k=0}^n s ( n , k ) B_k ( y - \lambda,
\lambda) ,
\]
where $s(k,l)$ is the Stirling number of the first
kind, that is,
$(-1)^{k-l} s(k,l)$ is the number of permutations of
$k$ elements which contain exactly $l$ permutation cycles,
$n\in\inte$; cf. Proposition \ref{p1} below.

The outline of this paper is as follows.
Section \ref{s2} contains our main results on
anticipative Girsanov identities and applications to
the Girsanov theorem.
In Section~\ref{exs} we consider some examples of anticipating
transformations to which this theorem can be applied;
this includes the adapted case as well as transformations that
act inside the convex hull generated by Poisson random measures,
given the positions of the extremal vertices.
In Section \ref{s2.1} we show that those results are consequences
of identities for multiple integrals and stochastic exponentials.
In Section \ref{s3} we review some results of \cite{prinv} (cf. also
\cite{priinvcr}) on the computation
of moments of Poisson stochastic integrals, and we derive
some of their corollaries to be applied in this paper.
In Section \ref{s4} we derive some combinatorial identities
that allow us, in particular, to rewrite the Charlier polynomials
into a form suitable to the use of moment identities.
Finally in Section \ref{s5} we prove the results of
Section \ref{s2.1}, and in Section \ref{s6} we make some remarks
on how the results of this paper can be connected to the Carleman--Fredholm
determinant.

\section{Main results}
\label{s2}
Let $\Omega^X$ denote the configuration space on
a $\sigma$-compact metric space $X$ with Borel
$\sigma$-algebra $\mathcal{B}(X)$, that is,
\[
\Omega^X = \bigl\{ \omega= ( x_i )_{i=1}^N \subset X, x_i\not= x_j \ \forall
i\not= j, N \in\inte\cup\{ \infty\} \bigr\}
\]
is the space of at most countable locally finite subsets
of $X$, endowed with the Poisson probability
measure $\pi_\sigma$ with $\sigma$-finite diffuse intensity
$\sigma(dx)$ on~$X$, which is characterized by its
Laplace transform
%
%
\begin{eqnarray}
\label{lp}
\psi_{\sigma} ( f ) &=& E_\sigma\biggl[ \exp\biggl( \int_X f(x) \bigl(
\omega(dx) - \sigma(dx)\bigr) \biggr) \biggr]\nonumber\\[-8pt]\\[-8pt]
&=& \exp\biggl( \int_X \bigl( e^{ f(x)} - f(x) -1\bigr)
\sigma(dx) \biggr) ,\nonumber
\end{eqnarray}
$f\in L^2_\sigma(X )$,
or by the Girsanov identity (\ref{1})
by taking $f(x) = \log( 1 + g (x ) )$,
$x \in X$, $g\in\mathcal{C}_c (X)$,
where $E_\sigma$ denotes the expectation under $\pi_\sigma$,
and $\mathcal{C}_c (X)$ is the space of continuous functions
with compact support in $X$.

Each element $\omega$ of $\Omega^X$ is identified to the
Radon point measure
\[
\omega= \sum_{i=1}^{\omega(X)} \epsilon_{x_i},
\]
where $\epsilon_x$ denotes the Dirac measure at $x\in X$,
and $\omega(X) \in\inte\cup\{ \infty\}$
denotes the cardinality of $\omega\in\Omega^X$.

Consider a measurable random transformation
\[
\tau\dvtx\Omega^X \times X \rightarrow X
\]
of $X$, let $\tau_* (\omega)$, $\omega\!\in\!\Omega^X$,
denote the image measure of $\omega(dx)$
by \mbox{$\tau(\omega, \cdot)\dvtx X\,{\to}\,X$}, that is,
%
%
\begin{equation}
\label{llk}
\tau_*\dvtx\Omega^X \to\Omega^X\vadjust{\goodbreak}
\end{equation}
maps
\[
\omega= \sum_{i=1}^{\omega(X)}
\epsilon_{x_i}\qquad
\mbox{to }
\tau_* (\omega) = \sum_{i=1}^{\omega(X)}
\epsilon_{\tau( \omega, x_i)}
.
\]
In other words, the random mapping $\tau_*\dvtx\Omega^X \to\Omega^X$
shifts each configuration point $x \in\omega$ according to
$x \mapsto\tau( \omega, x )$.

Let $D$ denote the finite difference gradient
defined on any random variable $F\dvtx\Omega^X \to\real$
as
\[
D_x F(\omega) = F(\omega\cup\{ x \} ) - F( \omega),\qquad
\omega\in\Omega^X,
x\in X,
\]
for any random variable $F\dvtx\Omega^X \to\real$; cf.
\cite{kree,yito,picard}.
The operator $D$ is continuous on the space
$\dee_{2,1}$ defined by the norm
\[
\Vert F \Vert_{2,1}^2 = \Vert F \Vert_{L^2 ( \Omega^X , \pi_\sigma)}^2
+ \Vert D F \Vert_{L^2( \Omega^X \times X , \pi_\sigma\otimes
\sigma)}^2,\qquad
F \in\dee_{2,1}.
\]
The next result is a Girsanov identity for random,
nonadapted shifts of Poisson configuration points,
obtained as a consequence of Proposition \ref{pc12.0}
below which is proved at the end of Section \ref{s2.1}.
Here we let $Y$ denote another metric space with Borel
$\sigma$-algebra $\mathcal{B} ( Y )$.
\begin{prop}
\label{cc11}
Assume that $\tau\dvtx\Omega^X \times X \to Y$ satisfies
the cyclic condition
%
%
\begin{equation}
\label{cyclic3}
D_{t_1} \tau( \omega, t_2 ) \cdots D_{t_k} \tau( \omega, t_1 ) =
0,\qquad
\sigma( dt_1), \ldots, \sigma( dt_k)\mbox{-a.e.},\qquad
\omega\in\Omega^X,\hspace*{-35pt}
\end{equation}
for all $k \geq2$, and let $g\dvtx Y \to\real$
be a measurable function such that
%
%
\begin{equation}
\label{st1}
E_\sigma\biggl[ e^{ \int_X | g ( \tau( \omega, x ) ) | \sigma(dx)
} \prod_{x\in\omega} \bigl( 1 + | g ( \tau( \omega, x ) ) | \bigr) \biggr] < \infty .
\end{equation}
Then we have
\[
E_\sigma\biggl[ e^{ - \int_X g( \tau( \omega, x ) ) \sigma(dx) }
\prod_{x\in\omega} \bigl( 1 + g ( \tau( \omega, x ) ) \bigr) \biggr] = 1 .
\]
\end{prop}

As a consequence of Proposition \ref{cc11},
if $\tau\dvtx\Omega^X \times X \to X$
satisfies (\ref{cyclic3}) and
$\tau( \omega, \cdot)\dvtx X \to Y$
maps $\sigma$ to a fixed measure $\mu$
on $(Y, \mathcal{B} (Y))$ for all
$\omega\in\Omega^X$, then we have
\begin{eqnarray*}
E_\sigma\biggl[ \prod_{x\in\omega} \bigl( 1 + g ( \tau( \omega, x ) ) \bigr) \biggr] & = &
e^{ \int_X g( \tau( \omega, x ) ) \sigma(dx) }
\\
& = & e^{ \int_Y g( y ) \mu( dy ) } ,\qquad g \in\mathcal{C}_c ( Y );
\end{eqnarray*}
hence $\tau_*\dvtx\Omega^X \to\Omega^X$
maps $\pi_\sigma$ to $\pi_\mu$, which recovers
Theorem 3.3 of \cite{prinv}.

Proposition \ref{cc11} then implies the following anticipating
Girsanov theorem, in which the Radon--Nikodym density is given
by a Dol\'{e}ans exponential.
\begin{corollary}
\label{cjds}
Assume that for all $\omega\in\Omega^X$,
$\tau( \omega,\cdot)\dvtx X \to X$ is invertible on $X$
and that for all $t_0,\ldots,t_k \in X$, $k \geq1$,
there exists $i\in\{ 0 ,\ldots,k\}$ such that
%
%
\begin{equation}
\label{cdfjkd}
D_{t_i} \tau( \omega, x ) = 0
\end{equation}
for all $x$ in a neighborhood of $t_{i+1 \mod k}$,
and that the density
\[
\phi( \omega, x ) : =
\frac{d\tau_*^{-1} ( \omega, \cdot) \sigma}{d\sigma} ( x ) - 1,\qquad
x \in X,
\]
exists for all $\omega\in\Omega^X$, with
%
%
\begin{equation}
\label{rmv}
E_\sigma\biggl[ e^{ ( 1 + \varepsilon) \int_X \phi( \omega, x )
\sigma(dx) } \prod_{x\in\omega} \bigl( 1 + \phi( \omega, x ) \bigr)^{1 +
\varepsilon} \biggr] < \infty
\end{equation}
for some $\varepsilon> 0$.
Then
we have the Girsanov identity
\[
E_\sigma\biggl[ F ( \tau_* ( \omega) ) e^{ - \int_X \phi( \omega, x )
\sigma(dx) } \prod_{x\in\omega} \bigl( 1 + \phi( \omega, x ) \bigr) \biggr] = E_\sigma[
F ]
\]
for all $F\in L^1(\Omega^X )$.
\end{corollary}
\begin{pf}
First we note that from (\ref{cdfjkd}), for all $\omega\in\Omega^X$
and $t_0,\ldots,t_k \in X$, $k \geq1$, there exists
$i\in\{0 ,\ldots,k\}$ such that
%
%
\begin{equation}
\label{jfklds}
D_{t_i} \tau( \omega, t_{i+1 \mod k} ) =
D_{t_i} \phi( \omega, t_{i+1 \mod k} ) = 0.
\end{equation}
Next from Proposition \ref{cc11}, for all $f\in\mathcal{C}_c (X)$ we have
\begin{eqnarray*}
&&
E_\sigma\biggl[ e^{ - \int_X f ( x ) \sigma(dx) - \int_X \phi( \omega, x )
\sigma(dx) } \prod_{x\in\omega} \bigl( 1 + f ( \tau( \omega, x ) ) \bigr) \bigl( 1 +
\phi( \omega, x ) \bigr) \biggr]
\\
&&\qquad =
E_\sigma\biggl[ e^{- \int_X f ( \tau( \omega, x ) ) ( 1 + \phi( \omega,
x ) ) \sigma(dx) - \int_X \phi( \omega, x ) \sigma(dx)}
\\
&&\qquad\quad\hphantom{E_\sigma\biggl[}{}
\times\prod_{x\in\omega} \bigl( 1 + f ( \tau( \omega, x ) ) + \phi(
\omega, x ) + f ( \tau( \omega, x ) ) \phi( \omega, x ) \bigr) \biggr]
\\
&&\qquad = 1
\end{eqnarray*}
by Proposition \ref{cc11}, since
\[
x \mapsto f ( \tau( \omega, x ) ) + \phi( \omega, x ) + f ( \tau(
\omega, x ) ) \phi( \omega, x )
\]
satisfies condition (\ref{cyclic3}) by (\ref{jfklds}).
We conclude by the density
in $L^1(\Omega^X )$ of linear combinations of
$F$ of the form
\[
F=\exp\biggl(- \int_X f ( x ) \sigma(dx)\biggr) \prod_{x\in\omega} \bigl( 1 + f ( x ) \bigr),
\qquad f \in\mathcal{C}_c (X).
\]
\upqed\end{pf}

Under the hypotheses of Corollary \ref{cjds},
if $\tau_*\dvtx\Omega^X \to\Omega^X$ is invertible
then the random transformation $\tau_*^{-1}\dvtx\Omega^X \to\Omega^X$
is absolutely continuous with respect to~$\pi_\sigma$,
with density
%
%
\begin{equation}
\label{rn}
\frac{d\tau^{-1}_* \pi_\sigma}{d\pi_\sigma} = e^{ - \int_X
\phi( \omega, x ) \sigma(dx) } \prod_{x\in\omega} \bigl( 1 + \phi( \omega, x
) \bigr) .
\end{equation}
In Corollary \ref{cjds}, condition (\ref{rmv}) actually requires
$\sigma( \tau( X ) )$ to be a.s. finite.

\section{Examples}
\label{exs}
In this section we present an example of a random nonadapted
transformation satisfying the hypotheses of Corollary \ref{cjds}.
First we note that condition (\ref{cyclic3}) is an extension of the usual
adaptedness condition, as it holds when $\tau\dvtx X \to X$ is adapted
to a given total binary relation $\preceq$ on~$X$.
Indeed, if $\tau\dvtx\Omega^X \times X \to X$ satisfies
\[
D_x \tau( \omega, y ) = 0,\qquad y \preceq x,
\]
then condition (\ref{cyclic3}) is satisfied
since for all $t_1,\ldots,t_k \in X$
there exists $i\in\{ 1,\ldots, k\}$ such that
$t_j \preceq t_i$, for all $1\leq j \leq k$;
hence $D_{t_i} \tau( \omega, t_j ) = 0$, \mbox{$1\leq j \leq k$}.
In this case, Corollary \ref{cjds} recovers a classical result
in the case where
$\tau\dvtx\allowbreak X \to X$ is deterministic or adapted;
cf., for example, Theorem 3.10.21 of~\cite{bichtelerbk}.

\begin{figure}[b]

\includegraphics{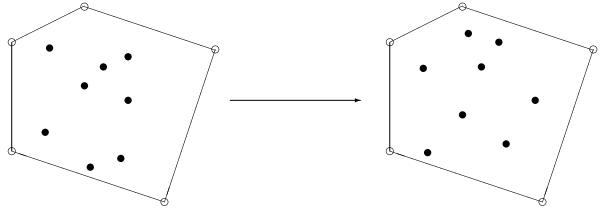}

\caption{Example of random transformation.}\label{fig1}
\end{figure}

Next, let
$X = \bar{B} (0,1)$ denote the closed unit ball in $\real^d$, with
$\sigma(dx)$ the Lebesgue measure.
For all $\omega\in\Omega^X$,
let $\mathcal{C} (\omega)$ denote the convex hull
of $\omega$ in~$X$ with interior
$\dot\mathcal{C} (\omega)$, and let
$\omega_e = \omega\cap
( \mathcal{C}(\omega) \setminus\dot\mathcal{C}(\omega) )$
denote the extremal vertices of $\mathcal{C} (\omega)$.
Consider a measurable\vspace*{1pt} mapping $\tau\dvtx\Omega^X \times X \to X$
such that for all $\omega\in\Omega^X$,
$\tau( \omega, \cdot)$ is measure preserving, maps
$\dot{\mathcal{C}}(\omega)$ to $\dot{\mathcal{C}}(\omega)$,
and for all $\omega\in\Omega^X$,
%
%
\begin{equation}
\label{fr1}
\tau( \omega, x ) = \cases{
\tau( \omega_e , x ),
&\quad$x\in\dot{\mathcal{C}}(\omega)$,\cr
x , &\quad$x\in X\setminus\dot{\mathcal{C}}(\omega)$,}
\end{equation}
that is, $\tau( \omega, \cdot)\dvtx X \to X$ modifies
only the inside points of the convex hull of $\omega$, depending on
the positions of its extremal vertices,
which are left invariant by $\tau( \omega, \cdot)$,
as illustrated in Figure \ref{fig1}.

%

Next, assume that $\tau( \omega, \cdot)\dvtx X \to X$
in (\ref{fr1}) has the form
\[
\tau( \omega, x ) = x + \psi( \omega_e , x ),\qquad
x \in X,
\]
for fixed $\omega\in\Omega^X$, where
$\psi(\omega_e , \cdot)\dvtx X \to X$ is a diffeomorphism such
that $\tau( \omega, \cdot)\dvtx X \to X$ is invertible for
all $\omega\in\Omega^X$; for example,
%
%
\begin{equation}
\label{psi}
\psi(\omega_e , x ) = u{\mathbf1}_{\mathcal{C} ( \omega)} ( x )
\frac{( d(x,\mathcal{C}(\omega) \setminus\dot{\mathcal{C}}(\omega)))^2}{
1 + ( d(x,\mathcal{C}(\omega) \setminus\dot{\mathcal{C}} (\omega ))
)^2},\qquad x \in X,
\end{equation}
with $u\in\real^d$ such that $\Vert u \Vert_d < 1/4$,
where $d(x,A)$ denotes the Euclidean distance from $x \in\real^d$
to the closed set $A\subset\real^d$.
Then the transformation $\tau\dvtx \Omega^X \times X \to X$
satisfies the hypotheses of
Corollary \ref{cjds} by Proposition \ref{pl} below, and
$\tau^* \dvtx \Omega^X \to\Omega^X$ is invertible with
\[
( \tau^* )^{-1} ( \omega) = \omega_e
\cup\bigcup_{x\in\omega\cap\dot\mathcal{C}( \omega)} \{ \tau^{-1} (
\omega_e , x ) \},\qquad \omega\in\Omega^X ;
\]
thus the associated Radon--Nikodym density (\ref{rn}) is given by taking
\[
\phi( \omega, x ) = \det\bigl( I_{\real^d} + \nabla_x \psi( \omega_e , x ) \bigr)
- 1,\qquad \omega\in X, x \in X.
\]
This quasi-invariance property is related to the intuitive fact that
a Poisson random measure remains Poisson
within its convex hull when its configuration points
are shifted given to the position of its extremal vertices; cf., for
example, \cite{davydov}.
\begin{prop}
\label{pl}
Assume that the random transformation $\tau\dvtx \Omega^X \times X \to X$
satisfies condition (\ref{fr1}).
Then $\tau$ satisfies the cyclic
condition (\ref{cdfjkd}) of Corollary \ref{cjds}.
\end{prop}
\begin{pf}
Let $t_1,\ldots,t_k\in X$.
First, if there exists $i \in\{1,\ldots,k\}$
such that $t_i \in\mathcal{C} ( \omega)$,
then for all $x\in X$ we have $t_i \in\mathcal{C} ( \omega\cup\{ x
\} )$,
and by Lemma \ref{ll00} below we get
\[
D_{t_i} \tau( \omega, x ) = 0,\qquad x \in X;
\]
thus (\ref{cdfjkd}) holds, and
we may assume that $t_i \notin\mathcal{C} ( \omega)$
for all $i = 1,\ldots,k$.
In this case, if
$t_{i+1 \mod k} \notin\mathcal{C} ( \omega\cup\{ t_i \} )$
for some $i=1,\ldots,k$, then by Lemma \ref{ll00} we have
\[
D_{t_i} \tau( \omega, t_{i+1 \mod k} ) = 0
;
\]
hence (\ref{cdfjkd}) holds since the set
$\mathcal{C} ( \omega\cup\{ t_i \} )$ is closed.
Next, if
$t_1 \in\mathcal{C} ( \omega\cup\{ t_k \} )$,
$t_k \in\mathcal{C} ( \omega\cup\{ t_{k-1} \} ),
\ldots,t_2 \in\mathcal{C} ( \omega\cup\{ t_1 \} )$,
then
we have
$t_1 \in\mathcal{C} ( \omega\cup\{ t_k \} )$ and
$t_k \in\mathcal{C} ( \omega\cup\{ t_1 \} )$,
which implies $t_1=t_k \notin\mathcal{C} ( \omega)$,
and we check that $D_{t_k} \tau( \omega, t_1 ) = 0$.
\end{pf}

Next we state and prove Lemma \ref{ll00} which has been used
above.
\begin{lemma}[\cite{prinv}]
\label{ll00}
For all $x, y \in X$ and $\omega\in\Omega^X$
we have
%
%
\begin{equation}
\label{f1}
x \in\mathcal{C} ( \omega\cup\{ y \} ) \quad\Longrightarrow\quad D_x
\tau( \omega, y ) = 0\vadjust{\goodbreak}
\end{equation}
and
%
%
\begin{equation}
\label{f2}
y \in\mathcal{C} ( \omega\cup\{ x \} ) \quad\Longrightarrow\quad D_x
\tau( \omega, y ) = 0.
\end{equation}
\end{lemma}
\begin{pf}
Let $x, y \in X$ and $\omega\in\Omega^X$.
First, if $y \notin\mathcal{C} ( \omega\cup\{ x \} )$
we have
$
\tau( \omega\cup\{x \} , y ) =
\tau( \omega, y ) = y $.
Next, if $x \in\mathcal{C} ( \omega\cup\{ y \} )$,
we can distinguish two cases:
\begin{longlist}
\item[(a)]
$x \in\mathcal{C}(\omega)$. In this case we have
$\mathcal{C}(\omega\cup\{ x \}) = \mathcal{C}(\omega) $;
hence
$
\tau( \omega\cup\{x \} , y ) =
\tau( \omega, y )
$
for all $y \in X$.
\item[(b)]
$x \in\mathcal{C} ( \omega\cup\{ y \} ) \setminus\mathcal{C} (
\omega)$.
If $y \in\mathcal{C} ( \omega\cup\{ x \} )$, then
$x=y \notin\dot\mathcal{C} ( \omega\cup\{x\} )$;
hence
$
\tau( \omega\cup\{x \} , y ) =
\tau( \omega, y )
$.
On the other hand if
$y \notin\mathcal{C} ( \omega\cup\{ x \} )$,
then
$
\tau( \omega\cup\{x \} , y ) =
\tau( \omega, y ) = y
$ as above.
\end{longlist}
We conclude that $D_x \tau( \omega, y ) = 0$
in both cases.
\end{pf}


\section{Multiple integrals and stochastic exponentials}
\label{s2.1}
$\!\!$The proofs of the above results will use properties of stochastic
exponentials and multiple stochastic integrals which are introduced
and proved in this section.
Let now
\[
I_n (f_n) (\omega) = \int_{\Delta_n} f_n(x_1,\ldots,x_n)
\bigl(\omega(dx_1)-\sigma(dx_1)\bigr) \cdots\bigl(\omega(dx_n)-\sigma(dx_n)\bigr)
\]
denote the multiple Poisson stochastic integral of
the symmetric function $f_n \in L^2_\sigma( X^n)$,
where
\[
\Delta_n = {\{(x_1,\ldots,x_n)\in X^n \dvtx x_i\not= x_j , \forall i\not=
j\}} ,
\]
with
\[
e^{ - \int_X g ( x ) \sigma(dx) } \prod_{x\in\omega} \bigl( 1 + g ( x ) \bigr) =
\sum_{n=0}^\infty\frac{1}{n!} I_n ( g^{\otimes n} )
\]
for $g\in L^2_\sigma( X )$ with bounded support,
where ``$\otimes$'' denotes the tensor product
of functions in $L^2_\sigma( X )$.
For all (possibly random) disjoint subsets $A_1,\ldots,A_n$
of $X$ with finite measure, we have the relation
%
%
\begin{equation}
\label{rel2}
I_N ( {\mathbf1}_{A_1^{k_1}} \circ \cdots \circ
{\mathbf1}_{A_n^{k_n}} ) = \prod_{i=1}^n C_{k_i} ( \omega( A_i ) , \sigma(
A_i ) )
\end{equation}
between the multiple Poisson integrals and the Charlier
polynomials,
where ``$\circ$'' denotes the symmetric tensor product
of functions in $L^2_\sigma( X )$ and $N=k_1+\cdots+ k_n$;
cf., for example, Proposition 6.2.9 in \cite{privaultbk2}.

Proposition \ref{cc11} will be proved using the following
Proposition \ref{pc12.0}
which is a restatement of Corollary \ref{pc12} below.
It provides a formula for the expectation of a multiple
stochastic integral of a time-changed function.
\begin{prop}
\label{pc12.0}
Assume that
$\tau\dvtx \Omega^X \times X \to Y$ satisfies
%
%
\begin{equation}
\label{fdsjklfds}
D_t \tau( \omega, t ) = 0,\qquad \omega\in\Omega^X, t \in
X.\vadjust{\goodbreak}
\end{equation}
Then for all symmetric step functions $g\dvtx Y^N \to\real$ of the
form
\[
g = \mathop{\sum_{k_1+\cdots+k_n = N}}_{1 \leq n \leq N} c_{k_1,\ldots,k_n}
{\mathbf1}_{B_{1,k_1}}^{\otimes^{k_1}} \circ \cdots \circ
{\mathbf1}_{B_{n,k_n}}^{\otimes k_n},
\]
where $N\geq1$ and $B_{1,k_1},\ldots,B_{n,k_n}$
are deterministic disjoint Borel subsets of~$Y$
and $c_{k_1,\ldots,k_n} \in\real$, we have
\begin{eqnarray*}
&&E_\sigma[ I_N ( {\mathbf1}_{A^N} ( \cdot) \tau^{\otimes N} ( \omega, \cdot)
) ] \\
&&\qquad= E_\sigma\biggl[ \int_{A^N} D_{t_1} \cdots D_{t_N} g ( \tau( \omega, t_1
) , \ldots, \tau( \omega, t_N ) ) \sigma( dt_1 ) \cdots\sigma( dt_N )
\biggr]
\end{eqnarray*}
for all compact subset $A\in\mathcal{B}(X)$ of $X$.
\end{prop}
\begin{pf}
It suffices to prove that for all deterministic disjoint Borel
subsets $B_1,\ldots,B_n$ of $Y$ we have
\begin{eqnarray*}
&&
E_\sigma\bigl[ I_N \bigl( {\mathbf1}_{A^N} {\mathbf1}_{\tau^{-1} ( B_1)}^{\otimes k_1}
\circ \cdots \circ {\mathbf1}_{\tau^{-1} ( B_n)}^{\otimes k_n} \bigr) \bigr]\\
&&\qquad =
E_\sigma\bigl[ I_N \bigl( {\mathbf1}_{A \cap\tau^{-1} ( B_1)}^{\otimes k_1} \circ
\cdots \circ {\mathbf1}_{A \cap\tau^{-1} ( B_n)}^{\otimes k_n}\bigr) \bigr]
\\
&&\qquad = E_\sigma\Biggl[ \int_{A^N} \Biggl( \prod_{i=1}^N D_{t_i} \Biggr) \bigl(
{\mathbf1}_{B_1^{k_1}} \otimes \cdots \otimes {\mathbf1}_{B_n^{k_n}} ( \tau(
\omega, t_1 ) , \ldots, \tau( \omega, t_N ) ) \bigr) \\
&&\qquad\quad\hspace*{183.3pt}{}\times\sigma( dt_1 )
\cdots\sigma( dt_N ) \Biggr]
\end{eqnarray*}
with $N = k_1 + \cdots+ k_n$,
and this is a direct consequence of relation (\ref{rel2})
above and Corollary \ref{pc12} below
applied to the random sets $A \cap\tau^{-1} ( B_1),
\ldots, A \cap\tau^{-1} ( B_n)$.
\end{pf}

As a particular case of Proposition \ref{pc12.0},
for $g = {\mathbf1}_B$ and $B \in\mathcal{B} (Y)$ such that
$\tau^{-1} ( B ) \subset A$ a.s., where $A$ is a fixed
compact subset of $X$, we have
%
%
\begin{eqnarray}
\label{adfgsdg}
&&
E_\sigma[ C_n ( \tau_* \omega( B ) , \tau_* \sigma( B ))]\nonumber\\[-8pt]\\[-8pt]
&&\qquad = E_\sigma\Biggl[ \int_{A^n} D_{s_1} \cdots D_{s_n} \prod_{p=1}^n
{\mathbf1}_B ( \tau( \omega, s_p ) ) \sigma( ds_1 ) \cdots\sigma( ds_n )
\Biggr],\nonumber
\end{eqnarray}
under condition (\ref{fdsjklfds}).
When $D_s {\mathbf1}_B ( \tau( \omega, t) ) $ is quasi-nilpotent in the
sense of condition (\ref{cyclic3}) above
for all $k \geq2$, $\omega\in\Omega^X$,
relation (\ref{adfgsdg}) and Lemma \ref{l12} below show that
\[
E_\sigma[ C_n ( \tau_* \omega( B ) , \tau_* \sigma( B ) ) ] = 0 ,
\]
and this extends (\ref{1ad}) as a particular case since
when $X=\real_+$, condition~(\ref{cyclic3}) holds
in particular when either
\[
D_s g ( \tau( \omega, t ) ) =0,\qquad 0 \leq s \leq t,
\]
or
\[
D_t g ( \tau( \omega, s ) ) =0,\qquad 0 \leq s \leq t,
\]
that is, when the process $\tau( \omega, t )$ is forward
or backward adapted with respect to the filtration generated
by the standard Poisson process $(N_t)_{t\in[0,T]}$.
\begin{pf*}{Proof of Proposition \ref{cc11}}
We take $g\dvtx Y \to\real$ to be the step function
\[
g ( t ) =\sum_{i=1}^m c_i {\mathbf1}_{B_i} ( t ),\qquad t \in Y,
\]
where $c_1,\ldots,c_m \in\real$ and
$B_1,\ldots, B_m \in\mathcal{B} (Y)$ are
disjoint Borel subsets of $Y$.
Then the expression
\[
C_n ( x , \lambda) = \sum_{k=0}^n x^k \sum_{l=0}^k \pmatrix{n \cr l} (
- \lambda)^{n-l} s(k,l),\qquad x,\lambda\in\real,
\]
for the Charlier polynomial of order $n \in\inte$,
shows that
\[
| C_n ( x , \lambda) | \leq \sum_{k=0}^n x^k \sum_{l=0}^k \pmatrix{n
\cr l} \lambda^{n-l} s(k,l) =C_n ( x , - \lambda),\qquad x , \lambda\geq0;
\]
hence
\[
\sum_{n=0}^\infty \frac{|r|^n}{n!} | C_n ( x , \lambda) | \leq e^{ | r
\lambda| } ( 1 + | r | )^x ,\qquad r \in\real,
\]
and letting $A\in\mathcal{B}(X)$ be a compact subset of $X$
we have
%
%
\begin{eqnarray}\label{bnd1}
&&E_\sigma\Biggl[ \sum_{n=0}^\infty\frac{1}{n!} | I_n ( {\mathbf1}_{A^n} ( \cdot)
g^{\otimes n} ( \tau^{\otimes n} ( \omega, \cdot) ) ) | \Biggr]
\nonumber\hspace*{-20pt}\\
&&\qquad = E_\sigma\Biggl[ \sum_{N=0}^\infty \Biggl| \mathop{\sum_{k_1+\cdots+k_n =
N}}_{n\geq0} \Biggl( \prod_{l=1}^n \frac{c_i^{k_i}}{k_i!} \Biggr) I_N \bigl(
{\mathbf1}_{A\cap\tau^{-1} ( B_1)}^{\otimes k_1} \circ \cdots \circ
{\mathbf1}_{A\cap\tau^{-1} ( B_n )}^{\otimes k_n} \bigr) \Biggr| \Biggr]
\nonumber\hspace*{-20pt}\\
&&\qquad = E_\sigma\Biggl[ \sum_{N=0}^\infty \Biggl| \mathop{\sum_{k_1+\cdots+k_n =
N}}_{n\geq0} \prod_{i=1}^n \frac{c_i^{k_i}}{k_i!} C_{k_i} \bigl( \omega\bigl(
A\cap\tau^{-1} ( B_i ) \bigr) , \sigma\bigl( A\cap\tau ^{-1} ( B_i ) \bigr) \bigr)\Biggr|\Biggr]
\nonumber\hspace*{-20pt}\\
&&\qquad \leq E_\sigma\Biggl[ \sum_{N=0}^\infty \mathop{\sum_{k_1+\cdots+k_n =
N}}_{n\geq0} \prod_{i=1}^n \frac{|c_i|^{k_i}}{k_i!} \bigl| C_{k_i} \bigl( \omega\bigl(
A\cap\tau^{-1} ( B_i ) \bigr) , \sigma\bigl( A\cap\tau ^{-1} ( B_i ) \bigr) \bigr) \bigr| \Biggr]
\nonumber\hspace*{-20pt}\\
&&\qquad \leq E_\sigma\Biggl[ \sum_{N=0}^\infty \mathop{\sum_{k_1+\cdots+k_n =
N}}_{n\geq0} \prod_{i=1}^n \frac{|c_i|^{k_i}}{k_i!} C_{k_i} \bigl( \omega\bigl(
A\cap\tau^{-1} ( B_i ) \bigr) , - \sigma\bigl( A\cap\tau^{-1} ( B_i ) \bigr) \bigr) \Biggr]
\nonumber\hspace*{-20pt}\\
&&\qquad = E_\sigma\Biggl[ \prod_{i=1}^n \sum_{k_i=0}^\infty
\frac{|c_i|^{k_i}}{k_i!} C_{k_i} \bigl( \omega\bigl( A\cap\tau^{-1} ( B_i ) \bigr) , -
\sigma\bigl( A\cap\tau^{-1} ( B_i ) \bigr) \bigr) \Biggr]\nonumber\hspace*{-20pt}\\
&&\qquad = E_\sigma\Biggl[ \prod_{i=1}^n \exp\bigl( |c_i| \sigma\bigl( A\cap\tau^{-1} ( B_i
) \bigr) \bigr) ( 1 + |c_i| )^{\omega( A\cap\tau^{-1} ( B_i ) ) } \Biggr]\nonumber\\
&&\qquad \leq E_\sigma\Biggl[ \prod_{i=1}^n \exp( |c_i| \tau_* \sigma( B_i ) ) ( 1
+ |c_i| )^{\tau_* \omega( B_i )} \Biggr]\hspace*{-20pt}\\
&&\qquad = E_\sigma\Biggl[ e^{ \int_X |g ( \tau( \omega, x ) ) | \sigma(dx) }
\prod_{x\in\omega} \bigl( 1 + |g ( \tau( \omega, x ) ) | \bigr) \Biggr]\nonumber\hspace*{-20pt}\\
&&\qquad < \infty.\nonumber
\end{eqnarray}
Consequently we can apply the Fubini theorem,
which shows that
\begin{eqnarray*}
&&E_\sigma\biggl[ e^{ - \int_A g( \tau( \omega, x ) ) \sigma(dx) } \prod_{x\in
A \cap\omega} \bigl( 1 + g ( \tau( \omega, x ) ) \bigr) \biggr] \\
&&\qquad= E_\sigma\Biggl[
\sum_{n=0}^\infty\frac{1}{n!} I_n ( {\mathbf1}_{A^n} ( \cdot) g^{\otimes n}
( \tau^{\otimes n} ( \omega, \cdot) ) ) \Biggr]
\\
&&\qquad = \sum_{n=0}^\infty\frac{1}{n!} E_\sigma[ I_n ( {\mathbf1}_{A^n} (
\cdot) g^{\otimes n} ( \tau^{\otimes n} ( \omega, \cdot) ) ) ]
\\
&&\qquad = \sum_{n=0}^\infty\frac{1}{n!} E_\sigma\Biggl[ \int_{A^n} D_{s_1} \cdots
D_{s_n} \prod_{p=1}^n g ( \tau( \omega, s_p ) ) \sigma( ds_1 )
\cdots\sigma( ds_n ) \Biggr]
\\
&&\qquad = 0
\end{eqnarray*}
by Proposition \ref{pc12.0}, provided
%
%
\begin{equation}
\label{lc00000}
\int_{A^n} D_{s_1} \cdots D_{s_n} \prod_{p=1}^n g (
\tau( \omega, s_p ) ) \sigma( ds_1 ) \cdots\sigma( ds_n ) =0,\qquad n\geq1,
\end{equation}
$\pi_\sigma(d\omega)$-a.s., which holds by Lemma \ref{l12} below since
$D_s \tau( \omega, t )$ is quasi-nilpo\-tent in the sense of (\ref{cyclic3}).
The extension from $A$ to $X$, and then from
$g$, a step function, to a measurable
function satisfying (\ref{st1}),
can be done by dominated convergence using
bound (\ref{bnd1}) above.
\end{pf*}

The above results can also be summarized in the
following general statement which is also proved in Section \ref{s5}
by the same argument as in the proof of Proposition \ref{cc11}.
\begin{prop}
\label{t1.1}
Assume that $\tau\dvtx \Omega^X \times X \to Y$ satisfies
\[
D_t \tau( \omega, t ) = 0,\qquad \omega\in\Omega^X, t \in X.
\]
Then for all bounded measurable functions $g \dvtx Y \to\real$
satisfying (\ref{st1}) we have
%
%
\begin{eqnarray}\label{cva}
&&
E_\sigma\Biggl[ e^{ - \int_X g( \tau( \omega, x ) ) \sigma(dx) }
\prod_{x\in\omega} \bigl( 1 + g ( \tau( \omega, x ) ) \bigr) \Biggr]
\nonumber\hspace*{-20pt}\\[-8pt]\\[-8pt]
&&\qquad = \sum_{n=0}^\infty \frac{1}{n!} E_\sigma\Biggl[ \int_{X^n} D_{s_1}
\cdots D_{s_n} \prod_{p=1}^n g ( \tau( \omega, s_p ) ) \sigma( ds_1 )
\cdots\sigma( ds_n ) \Biggr] ,\nonumber\hspace*{-20pt}
\end{eqnarray}
provided
%
%
\begin{equation}
\label{bnd1.1}\qquad
\sum_{n=0}^\infty \frac{1}{n!} E_\sigma\Biggl[ \int_{X^n} \Biggl|
D_{s_1} \cdots D_{s_n} \prod_{p=1}^n g ( \tau( \omega, s_p ) ) \sigma(
ds_1 ) \cdots\sigma( ds_n ) \Biggr| \Biggr] < \infty.
\end{equation}
\end{prop}

In the next lemma we show that relation (\ref{lc00000})
is satisfied provided $D_s \tau( \omega, t)$ satisfies the
cyclic condition (\ref{cyclic2}) below.
\begin{lemma}
\label{l12}
Let $N\geq1$, and assume that
$\tau\dvtx \Omega^X \times X \to X$
satisfies the cyclic condition
%
%
\begin{equation}
\label{cyclic2}
D_{t_0} \tau( \omega, t_1 ) \cdots D_{t_k} \tau(
\omega, t_0 ) = 0,\qquad \omega\in\Omega^X, t_0,\ldots,t_k \in X,
\end{equation}
for $k = 1, \ldots, N$.
Then we have
\[
D_{t_0} \cdots D_{t_k} \prod_{p=0}^k g ( \tau( \omega, t_p ) ) =0,\qquad
t_0, \ldots, t_k \in X,
\]
for $k = 1, \ldots, N$.
\end{lemma}
\begin{pf}
We use the relation
%
%
\begin{eqnarray}
\label{defdlta}
&&
D_{s_0} \cdots D_{s_j} \prod_{p=0}^n g ( \tau( \omega,
s_p ) )\nonumber\\[-8pt]\\[-8pt]
&&\qquad = \sum_{ { \Theta_0 \cup\cdots\cup\Theta_n = \{ 0 , 1 , \ldots,
j \} } } D_{\Theta_0} g ( \tau( \omega, s_0 ) ) \cdots D_{\Theta_n} g (
\tau( \omega, s_n ) ) ,\nonumber
\end{eqnarray}
$s_0,\ldots,s_n \in X$,
where $ D_\Theta: = \prod_{j\in\Theta} D_{s_j}$ when
$\Theta\subset\{0,1,\ldots,j\}$, $0 \leq j \leq n$,
which follows from the product rule
%
%
\begin{equation}
\label{prdr}
D_t ( FG ) = FD_t G + G D_t F + D_t F D_t G,\qquad t\in X,
\end{equation}
which is satisfied by $D_t$ as a finite difference operator.
Without loss of generality we may assume that
$\Theta_0 \not= \varnothing, \ldots, \Theta_j \not= \varnothing$
and $\Theta_k \cap\Theta_l = \varnothing$, $0 \leq k \not=l \leq j$.
In this case we can construct a sequence $(k_1,\ldots,k_i)$ by
choosing
\[
0\not= k_1\in\Theta_0,\qquad k_2 \in\Theta_{k_1},
\ldots,
k_{i-1} \in\Theta_{k_{i-2}},
\]
until $k_i=0\in\Theta_{k_{i-1}}$ for some $i\in\{2,\ldots,j \}$
since $\Theta_0 \cap\cdots\cap\Theta_j =\varnothing$
and
$\Theta_0 \cup\cdots\cup\Theta_j = \{ 0 , 1 , \ldots, j \}$.
Hence by (\ref{cyclic2}) we have
\begin{eqnarray*}
&&D_{s_{k_1}} g ( \tau( \omega, s_{s_0} ) ) D_{s_{k_2}} g ( \tau( \omega,
s_{s_{k_1}} ) ) \cdots \\
&&\qquad{}\times D_{s_{k_{i-1}}} g ( \tau( \omega,
s_{s_{k_{i-2}}} ) ) D_{s_0} g ( \tau( \omega, s_{s_{k_{i-1}}} ) ) = 0
\end{eqnarray*}
by (\ref{cyclic2}), which implies
\begin{eqnarray*}
&&D_{\Theta_0} g ( \tau( \omega, s_0 ) ) D_{\Theta_{k_1}} g ( \tau(
\omega, s_{k_1} ) ) \cdots \\
&&\qquad{}\times D_{\Theta_{k_{i-2}}} g ( \tau( \omega,
s_{k_{i-2}} ) ) D_{\Theta_{k_{i-1}}} g ( \tau( \omega, s_{k_{i-1}} ) )
= 0 ,
\end{eqnarray*}
since
\[
( k_1 , \ldots, k_{i-1} , 0 ) \in \Theta_0 \times\Theta_{k_1}
\times\cdots\times\Theta_{k_{i-1}} .
\]
\upqed\end{pf}

\section{Moment identities for Poisson integrals}
\label{s3}

In this section we state some results obtained in
\cite{prinv} on the moments of Poisson stochastic
integrals, and we reformulate them in view of our
applications to Girsanov identities and to
random Charlier polynomial functionals.

The Poisson--Skorohod integral operator $\delta$ is defined
on any measurable
process $u \dvtx \Omega^X \times X \to\real$
by the expression
%
%
\begin{equation}
\label{ui}
\delta( u ) = \int_X u ( \omega\setminus\{ t \} , t ) \bigl(
\omega( dt ) - \sigma(dt) \bigr) ,
\end{equation}
provided
$E_\sigma[ \int_X | u ( \omega, t ) | \sigma(dt ) ] <
\infty$;
cf., for example, \cite{nualartvives,privaultbk2}.

Note that if $D_t u_t = 0$, $t\in X$, and
in particular when applying (\ref{ui})
to $u \in L^1_\sigma( X )$ a deterministic function, we have
%
%
\begin{equation}
\label{ui2}
\delta( u ) = \int_X u ( t ) \bigl( \omega( d t ) - \sigma(d t )\bigr),
\end{equation}
that is, $\delta( u )$
with the compensated Poisson--Stieltjes integral of $u$.
In addition, if $X = \real_+$ and $\sigma(dt) = \lambda_t \,dt$,
we have
%
%
\begin{equation}
\label{*9}
\delta( u ) = \int_0^\infty u_t (dN_t - \lambda_t \,dt)
\end{equation}
for\vspace*{1pt} all square-integrable predictable processes $(u_t)_{t\in\real_+}$,
where $N_t = \omega([0,t])$, \mbox{$t\in\real_+$}, is a Poisson process with
intensity $\lambda_t >0$; cf., for instance, the example on page 518 of
\cite{picard}.\vadjust{\goodbreak}

From Corollaries 1 and 5 in \cite{picard} or
Proposition 6.4.3 in \cite{privaultbk2}
the operators~$D$ and $\delta$ are closable
and satisfy the duality relation
%
%
\begin{equation}
\label{dr1}
E_\sigma\bigl[ \langle D F , u \rangle_{L^2_\sigma( X )}\bigr] =
E_\sigma[ F \delta( u ) ],
\end{equation}
which can be seen as a formulation of the Mecke
\cite{jmecke}
identity for Poisson random measures,
on their $L^2$ domains
$\Dom( \delta)
\subset L^2 ( \Omega^X \times X , \pi_\sigma\otimes\sigma)$
and
$\Dom( D ) = \dee_{2,1}
\subset L^2 ( \Omega^X , \pi_\sigma)$
under the Poisson measure $\pi_{\sigma}$ with intensity~$\sigma$.

The operator $\delta$ is continuous on
the space $\lee_{2,1} \subset\Dom( \delta)$
defined by the norm
\[
\Vert u \Vert_{2,1}^2 = E_\sigma \biggl[ \int_X | u_t |^2 \sigma(d t ) \biggr] +
E_\sigma \biggl[ \int_X | D_s u_t |^2 \sigma(d s ) \sigma(d t ) \biggr] ,
\]
%
and it satisfies the Skorohod isometry
%
%
\begin{equation}
\label{skois}\qquad
E_\sigma[ \delta( u )^2 ] = E_\sigma\biggl[ \int_X | u_t |^2
\sigma( dt ) \biggr] + E_\sigma\biggl[ \int_X \int_X D_s u_t D_t u_s \sigma( ds )
\sigma( dt ) \biggr]
\end{equation}
for any $u\in\lee_{2,1}$; cf. Corollary 4 and pages 517 and 518 of
\cite{picard}.

In addition, from (\ref{ui}), for any $u\in\Dom(\delta)$
we have the commutation relation
%
%
\begin{equation}
\label{commrel}
D_t \delta(u) = \delta(D_t u) + u_t,\qquad t \in X,
\end{equation}
or
%
%
\begin{equation}
\label{commrel2}
( I + D_t ) \delta(u) = \delta\bigl( ( I + D_t) u\bigr) + u_t,\qquad t
\in X,
\end{equation}
provided $D_t u \in\lee_{2,1}$, $t\in X$.

The following lemma relies on the application
of relations (\ref{dr1}) and (\ref{commrel}),
and extends (\ref{skois}) to powers of order
greater than two; cf. Lemma 2.4 in \cite{prinv}.
\begin{lemma}[\cite{prinv}]
\label{l5}
Let $u \in\lee_{2,1}$ be such that $D_t u \in\lee_{2,1}$,
$t\in X$,
$\delta( u )^n \in\dee_{2,1}$,
and
\begin{eqnarray*}
E_\sigma\biggl[ \int_X | u_t |^{n-k+1} \bigl| \delta\bigl( ( I + D_t ) u \bigr) \bigr|^k \sigma(
dt ) \biggr] &<& \infty,\\[-2.5pt]
E_\sigma\biggl[ | \delta( u ) |^k \int_X | u_t |^{n-k+1}
\sigma( dt ) \biggr] &<& \infty ,
\end{eqnarray*}
$0 \leq k \leq n$.
Then we have
\begin{eqnarray*}
E_\sigma[ \delta( u )^{n+1} ] & = & \sum_{k=0}^{n-1} \pmatrix{n \cr k}
E_\sigma\biggl[ \delta( u )^k \int_X u_t^{n-k+1} \sigma( dt ) \biggr]
\\[-2.5pt]
& &{} + \sum_{k=1}^n \pmatrix{n \cr k} E_\sigma\biggl[ \int_X u_t^{n-k+1} \bigl(
\delta\bigl( ( I + D_t ) u \bigr)^k - \delta( u )^k \bigr) \sigma( dt ) \biggr]
\end{eqnarray*}
for all $n\geq1$.
\end{lemma}

When $h$ is a deterministic
function, Lemma \ref{l5}
yields the recursive covariance identity
%
%
\begin{equation}
\label{rmi}
E_\sigma[ \delta( h )^{n+1} ] =\sum_{k=1}^n \pmatrix{n \cr
k} \int_X h^{k+1} ( t ) \sigma( dt ) E_\sigma[\delta( h )^{n-k}],\qquad
n\geq0,\hspace*{-25pt}
\end{equation}
for the Poisson stochastic integral
\[
\delta( h ) = \int_X h(x) \bigl( \omega(dx) - \sigma(dx) \bigr).
\]
By induction, (\ref{rmi}) shows that
the moments of the above Poisson stochastic integral can be
computed as
%
%
\begin{equation}
\label{eee}
E_\sigma [ \delta( h)^n ] = \sum_{a=1}^{n-1} \sum_{0 = k_1
\ll\cdots\ll k_{a+1} = n } \prod_{l=1}^a \pmatrix{k_{l+1} - 1 \cr k_l }
\prod_{l=1}^a \int_X h^{k_{l+1}-k_l} \,d \sigma\hspace*{-25pt}
\end{equation}
for all $n\geq1$ and deterministic
$ h \in\bigcap_{p=2\wedge n}^n L^p_\sigma(X)$,
where $a \ll b$ means $a<b-1$, $a,b\in\inte$.
This result can also be recovered from the relation
%
%
\begin{equation}
\label{kmn}
E_\sigma [ \delta( h)^n ] = \sum_{d=1}^n \sum_{B_1,\ldots,
B_d} \kappa_{|B_1|} \cdots\kappa_{|B_d|},
\end{equation}
where the sum runs over all partitions of $\{ 1 , \ldots, n \}$,
$|B_i|$ denotes the cardinality of~$B_i$,
and $\kappa_1 = 0$, $\kappa_n = \int_X h^n (t) \sigma(dt)$,
$n\geq2$, denote the cumulants of~$\delta(h)$.

In particular, relations (\ref{eee}) and (\ref{kmn}) yield the identity
%
%
\begin{equation}
\label{as0}
E_\lambda [ ( Z - \lambda)^n ] = \sum_{a=0}^n \lambda^a S_2
( n , a )
\end{equation}
for the central moments of a Poisson random variable $Z$
with intensity $\lambda$, where
\[
S_2 ( n , a ) : = \sum_{0 = k_1 \ll\cdots\ll k_{a+1} = n }
\prod_{l=1}^a \pmatrix{k_{l+1} - 1 \cr k_l } ,
\]
represents the number of partitions of a set of size
$m$ into $a$ subsets of size at least~$2$.

In the sequel we let
%
%
\begin{eqnarray}
\label{defin}
&&C (l_1,\ldots,l_a , b ) \nonumber\\[-8pt]\\[-8pt]
&&\quad= \sum_{0 = r_{b+1} < \cdots< r_0
= a+b+1} \prod_{q=0}^b \prod_{p=r_{q+1} - ( b - q-1 )}^{r_q-1-(b-q)}
\pmatrix{l_1 + \cdots+ l_p + q - 1 \cr l_1 + \cdots+ l_{p-1} + q}
,\nonumber\hspace*{-25pt}
\end{eqnarray}
which represents the number of partitions of a set
of $l_1+\cdots+ l_a+b$ elements into $a$
subsets of lengths $l_1,\ldots,l_a$ and $b$
singletons.
We will need the following result; cf. Theorem 5.1 of
\cite{prinv}.\vadjust{\goodbreak}
\begin{theorem}[\cite{prinv}]
\label{l22}
Let $F\dvtx\Omega^X \to\real$ be a bounded random variable, and
let $u\dvtx\Omega^X \times X \to\real$ be a bounded process
with compact support in $X$.
For all $n\geq0$ we have
\begin{eqnarray*}
&&E_\sigma[ F \delta_\sigma( u )^n ] \\[-3pt]
&&\qquad= \sum_{a=0}^n \sum_{b = 0}^{n-a}
(-1)^b \mathop{\sum_{l_1 + \cdots+ l_a = n - b }}_{ l_1,\ldots,l_a \geq1} C (
l_1 , \ldots, l_a , b ) \\[-3pt]
&&\qquad\quad\hspace*{29.3pt}{}\times E_\sigma\Biggl[ \int_{X^{a+b}} \Biggl( \prod_{i=1}^a ( I +
D_{s_i} ) F \Biggr)
\Biggl( \prod_{q=a+1}^{a+b} \prod_{i=1}^a ( I + D_{s_i} ) u_{s_q}\Biggr)\\[-3pt]
&&\hspace*{37.7pt}\hspace*{48pt}\qquad\quad{} \times
\prod_{p=1}^a \Biggl( \mathop{\prod_{i=1}}_{i\not= p}^a ( I + D_{s_i} ) u_{s_p}
\Biggr)^{l_p} \sigma( ds_1 ) \cdots\sigma( ds_{a+b} ) \Biggr]
.\vspace*{-3pt}
\end{eqnarray*}
\end{theorem}

In the above proposition, by saying that
$u\dvtx\Omega^X \times X \to\real$ has a compact support in $X$ we
mean that there exists a compact $K \in\mathcal{B}(X)$ such that
$u(\omega, x ) = 0$ for all $\omega\in\Omega^X$ and $x \in X
\setminus K$.

In particular when $u = {\mathbf1}_A$
is a (random) indicator function we get the~fol\-lowing proposition,
which will be used to prove Proposition \ref{p12}
below.
We~let
%
%
\begin{equation}
\label{snc}
S(n,c) = \frac{1}{c!} \sum_{l=0}^ c (-1)^{c-l} \pmatrix{c
\cr l} l ^n
\end{equation}
denote the Stirling number of the
second kind, that is, the number of ways to partition a set of $n$
objects
into $c$ nonempty subsets.
In the next proposition, which is an application of Theorem \ref{l22},
the random indicator function
$(x,\omega) \mapsto{\mathbf1}_{A(\omega)}(x)$ on $\Omega^X \times X$
denotes a measurable process
$u\dvtx\Omega^X \times X \to\real$ such that
$u^2(\omega, t) = u (\omega, t)$, $\omega\in\Omega^X$, $t\in X$.\vspace*{-3pt}

\begin{prop}
\label{l221}
Let $F\dvtx\Omega^X \to\real$ be a bounded random variable, and
consider a measurable random indicator function
$(x,\omega) \mapsto {\mathbf1}_{A(\omega)}(x)$ on $\Omega^X\times X$,
with compact support in $X$.
Then for all $n\geq0$ we have
\begin{eqnarray*}
&&E_\sigma[ F \delta( {\mathbf1}_A )^n ]\\[-3pt]
&&\qquad = \sum_{c=0}^n \sum_{a=0}^c (-1)^a
\pmatrix{n \cr a} S ( n - a , c - a )
\\[-3pt]
&&\hspace*{28.9pt}\qquad\quad{}\times E_\sigma\Biggl[ \int_{X^a} \Biggl( \prod_{i=1}^a ( I + D_{s_i} ) ( F \sigma( A
)^{c-a} ) \Biggr) \\[-3pt]
&&\hspace*{77.3pt}\qquad\quad{}\times\prod_{p=1}^a \mathop{\prod_{i=1}}_{i\not= p}^a ( I + D_{s_i} )
{\mathbf1}_A ( s_p ) \sigma( ds_1 ) \cdots\sigma( ds_a ) \Biggr]
.\vspace*{-3pt}\vadjust{\goodbreak}
\end{eqnarray*}
\end{prop}
\begin{pf}
Taking $u = {\mathbf1}_A$ in Theorem \ref{l22} yields
\begin{eqnarray*}
&&
E_\sigma[ F ( \delta( u ) )^n ] \\[2pt]
&&\qquad= \sum_{a=0}^n \sum_{b = 0}^{n-a}
(-1)^b \mathop{\sum_{l_1 + \cdots+ l_a = n - b}}_{l_1,\ldots,l_a \geq1} C (
l_1 , \ldots, l_a , b ) \\[2pt]
&&\hspace*{29.3pt}\qquad\quad{}\times E_\sigma\Biggl[ \int_{X^{a+b}} \Biggl( \prod_{i=1}^a ( I +
D_{s_i} ) F \Biggr)
\\[2pt]
&&\hspace*{57pt}\hspace*{29.3pt}\qquad\quad{}\times
\prod_{p=1}^{a+b} \mathop{\prod_{i=1}}_{i\not= p}^a ( I + D_{s_i} )
{\mathbf1}_A ( s_p ) \sigma( ds_1 ) \cdots\sigma( ds_{a+b} ) \Biggr]
\\[2pt]
&&\qquad = \sum_{c=0}^n \sum_{a=0}^c (-1)^a \pmatrix{n \cr a} S ( n - a , c
- a )
\\[2pt]
&&\hspace*{29.3pt}\qquad\quad{}\times E_\sigma\Biggl[ \int_{X^c} \Biggl( \prod_{i=1}^a ( I + D_{s_i} ) F
\Biggr)\\[2pt]
&&\hspace*{76.7pt}\qquad\quad{}\times
\prod_{p=1}^c \mathop{\prod_{i=1}}_{i\not= p}^a ( I + D_{s_i} ) {\mathbf1}_A (
s_p ) \sigma( ds_1 ) \cdots\sigma( ds_c ) \Biggr]
\\[2pt]
&&\qquad = \sum_{c=0}^n \sum_{a=0}^c (-1)^a \pmatrix{n \cr a} S ( n - a , c
- a )\\[2pt]
&&\hspace*{29.3pt}\qquad\quad{}\times E_\sigma\Biggl[ \int_{X^a} \Biggl( \prod_{i=1}^a ( I + D_{s_i} ) ( F \sigma( A
)^{c-a} ) \Biggr) \\[2pt]
&&\qquad\quad\hspace*{76.7pt}{}\times\prod_{p=1}^a \mathop{\prod_{i=1}}_{i\not= p}^a ( I + D_{s_i} )
{\mathbf1}_A ( s_p ) \sigma( ds_1 ) \cdots\sigma( ds_a ) \Biggr] ,
\end{eqnarray*}
after checking that we have
\[
\pmatrix{n \cr b} S ( n - b , a ) = \mathop{\sum_{l_1 + \cdots+ l_a = n - b
}}_{l_1,\ldots,l_a \geq1} C ( l_1 , \ldots, l_a , b ) ,
\]
which is the number of partitions of
a set of $n$ elements into $a$ nonempty subsets
and one subset of size $b$.
\end{pf}

When the set $A$ is deterministic, Proposition \ref{l221} yields
\[
E_\lambda[ ( Z - \lambda)^n ] = \sum_{c=0}^n \lambda^c \sum_{a=0}^c
(-1)^a \pmatrix{n \cr a} S ( n - a , c - a )
\]
for the central moments of a Poisson random variable
$Z = \omega( A ) $ with intensity $\lambda= \sigma(A)$,
which, from (\ref{as0}), shows the combinatorial identity
%
%
\begin{equation}
\label{id}
S_2 ( n , c ) = \sum_{a = 0}^c (-1)^a \pmatrix{n \cr a } S (
n - a , c - a ) .
\end{equation}

\section{Poisson moments and polynomials}
\label{s4}
As mentioned in the \hyperref[sec1]{Introduction}
we need to introduce another family of polynomials
whose generating function and associated combinatorics will be better
adapted to our approach, %
making it possible to apply the moment identities of Proposition \ref{l221}
and the integration by parts formula (\ref{dr1}).

In terms of polynomials
the identity (\ref{adfgsdg}) is easy to check for $n=1$ and $n=2$,
in which case we have
\[
C_1 ( \omega(A ), \sigma(A) ) = \omega( A) - \sigma( A ) = \delta(
{\mathbf1}_A )
\]
and
%
%
\begin{eqnarray}\label{c1}
C_2 ( \omega( A ) , \sigma( A ) )
&=& \bigl( \omega( A ) - \sigma( A ) \bigr)^2
- \bigl( \omega( A ) - \sigma( A ) \bigr) - \sigma( A )
\nonumber\\[-8pt]\\[-8pt]
&=& \delta( {\mathbf1}_A )^2 - \delta( {\mathbf1}_A ) - \sigma( A )
,\nonumber
\end{eqnarray}
hence
\begin{eqnarray*}
E_\sigma[ C_2 ( \omega( A ) , \sigma( A ) ) ]
&=& E_\sigma[ \delta(
{\mathbf1}_A )^2 ] - \sigma( A )
\\
&=& E_\sigma\biggl[ \int_X \int_X D_s {\mathbf1}_A ( t ) D_t {\mathbf1}_A ( s )
\sigma( ds ) \sigma( dt ) \biggr]
\end{eqnarray*}
from the Skorohod isometry (\ref{skois}).

In the sequel we will need to extend the above calculations and
the proof of~(\ref{adfgsdg}) to Charlier polynomials
$C_n ( x , \lambda)$ of all orders.
For this, in Section~\ref{s5} we will use the moment identities for the
Skorohod integral $\delta( {\mathbf 1}_A )$ of Proposition~\ref{l221},
and for this reason we will need to rewrite
$C_n ( \omega( A ) , \sigma( A ) )$, a linear combination of
polynomials of the form $B_n ( \delta( {\mathbf1}_A ) , \sigma( A ) )$,
where $B_n ( x , \lambda)$ is another polynomial of degree $n$.
This construction is done
using Stirling numbers and combinatorial arguments;
cf. Proposition \ref{p1} below.

In other words, instead of using the identity (\ref{1})
we need its Laplace form (\ref{lp}), that is,
%
%
\begin{equation}
\label{l012}
E_\sigma\biggl[ \exp \biggl( \delta( f ) - \int_X \bigl( e^{f(x)} - f(x)
-1\bigr) \sigma(dx) \biggr) \biggr] = 1 ,
\end{equation}
obtained from (\ref{1}) by taking
\[
f(x) = \log\bigl( 1 + g (x ) \bigr),\qquad
x \in X.\vadjust{\goodbreak}
\]
In particular when $f = {\mathbf1}_A$ with $A\in\mathcal{B} (X)$
a fixed compact subset of $X$, relation (\ref{l012}) reads
%
%
\begin{equation}
\label{rw}
E_\sigma\bigl[e^{t \delta( {\mathbf1}_A ) - \sigma(A) (e^t - t - 1)}\bigr]
=1,\qquad t \in\real,
\end{equation}
where $\delta( {\mathbf1}_A ) = \omega( A ) - \sigma( A)$
is a compensated Poisson random variable with intensity
$\sigma(A)>0$.

We let $(B_n (x , \lambda))_{n\in\inte}$ denote the family of
polynomials defined by the generating function
%
%
\begin{equation}
\label{gf}
e^{ t y - \lambda(e^t - t - 1 ) } =
\sum_{n=0}^\infty\frac{t^n}{n!} B_n ( y , \lambda) ,\qquad t \in\real,
\end{equation}
for all $y,\lambda\in\real$.
This definition implies in particular that
\[
B_n ( y , - \lambda) = E_\lambda[ ( Z + y - \lambda)^n],
\]
where $Z$ is a Poisson random variable with intensity $\lambda>0$,
and
%
%
\begin{equation}
\label{ct}
B_n ( y , \lambda) =\sum_{k=0}^n \pmatrix{n \cr k} y^k
B_{n-k} ( 0 , \lambda),\qquad \lambda\in\real, n \in\inte.
\end{equation}
For example, one has that $B_1 ( y , \lambda) = y$ and
$B_2 ( y , \lambda)
=
y^2
-
\lambda$;
hence (\ref{c1}) reads
\[
C_2 ( x , \lambda) = B_2 ( x - \lambda, \lambda) - B_1 ( x - \lambda,
\lambda) ,
\]
and these relations will extended to all polynomial degrees
in Proposition~\ref{p1} below.

In addition, the definition of $B_n ( x , \lambda)$
generalizes that of the Bell (or Touchard)
polynomials $B_n(\lambda)$ defined by the generating function
\[
e^{ \lambda(e^t-1)} = \sum_{n=0}^\infty \frac{t^n}{n!} B_n ( \lambda) ,
\]
which satisfy
%
%
\begin{equation}
\label{ws}
B_n ( \lambda) = B_n ( \lambda, - \lambda ) = E_\lambda[ Z^n
] = \sum_{c=0}^n \lambda^c S(n,c) ,
\end{equation}
where $Z$ is a Poisson random variable with intensity $\lambda>0$;
cf., for example, Proposition 2 of \cite{boyadzhiev} or Section 3.1 of
\cite{dinardo}.

Next we show that the Charlier polynomials
$C_n ( x , \lambda)$ with exponential generating function
\[
e^{ - \lambda t } ( 1 + t )^x = \sum_{n=0}^\infty\frac{t^n}{n!} C_n ( x
, \lambda) ,\qquad x,t, \lambda\in\real,
\]
are dual to the generalized Bell polynomials
$B_n ( x - \lambda, \lambda)$ under the Stirling transform.
\begin{prop}
\label{p1}
We have the relations
\[
C_n ( y , \lambda) = \sum_{k=0}^n s ( n , k ) B_k ( y - \lambda,
\lambda)
\]
and
\[
B_n ( y , \lambda) = \sum_{k=0}^n S ( n , k ) C_k (
y + \lambda, \lambda) ,
\]
$y,\lambda\in\real$, $n\in\inte$.
\end{prop}
\begin{pf}
For the first relation, for all fixed $y, \lambda\in\real$ we let
\[
A ( t ) = e^{ - \lambda t } ( 1 + t )^{y+\lambda} =
\sum_{n=0}^\infty\frac{t^n}{n!} C_n ( y+\lambda, \lambda) ,\qquad t \in\real,
\]
and note that
\[
A ( e^t - 1 ) = e^{ t ( y+\lambda) - \lambda( e^t - 1 ) } =
\sum_{n=0}^\infty\frac{t^n}{n!} B_n ( y , \lambda) ,\qquad t \in\real,
\]
which implies
\[
B_n ( y , \lambda) = \sum_{k=0}^n S ( n , k ) C_k ( y + \lambda,
\lambda),\qquad
n \in\inte,
\]
(see, e.g., \cite{sloane}, page 2).
The second part can be proved
by inversion using Stirling numbers of the first kind,
as
\begin{eqnarray*}
\sum_{k=0}^n S ( n , k ) C_k ( y + \lambda, \lambda) & = & \sum_{k=0}^n
\sum_{l=0}^k S ( n , k ) s ( k , l ) B_l ( y , \lambda)
\\
& = & \sum_{l=0}^n B_l ( y , \lambda) \sum_{k=l}^n S ( n , k ) s ( k ,
l )
\\
& = & B_n ( y , \lambda)
\end{eqnarray*}
from the inversion formula
%
%
\begin{equation}
\label{inv}
\sum_{k=l}^n S ( n , k ) s ( k , l ) = {\mathbf1}_{\{ n = l \}}
,\qquad n,l \in\inte,
\end{equation}
for Stirling numbers; cf., for example, page 825 of \cite{Hand1972}.
\end{pf}

The combinatorial identity proved in the next lemma will
be used in Section \ref{s5} for the proof of Proposition \ref{p12}.
For $b=0$ it yields the
identity
%
%
\begin{equation}
\label{fd}
S ( n , a ) = \sum_{c=0}^a \pmatrix{n \cr c} S_2 ( n - c , a
- c ) ,
\end{equation}
which is the inversion formula of (\ref{id}),
and has a natural interpretation by stating that
$S_2 ( m , b )$
is the number of partitions of a set of $m$ elements
made of $b$ sets of cardinal greater or equal to $2$.
\begin{lemma}
\label{ll}
For all $a,b \in\inte$ we have
\[
\pmatrix{a+b \cr a} S ( n , a + b ) = \sum_{l=0}^b \sum_{k=l}^n
\pmatrix{n \cr k} \pmatrix{k \cr l} S ( k - l , a ) S_2 ( n-k , b-l ) .
\]
\end{lemma}
\begin{pf}
This identity can be proved by a combinatorial argument.
For each value of $k = 0 , \ldots, n$ one chooses
a subset of $\{ 1 , \ldots, n \}$ of size $k - l$
which is partitioned into $a$ nonempty subsets,
the remaining set of size $n+l-k$ being partitioned
into $l$ singletons and $b-l$ subsets of size at
least $2$.
In this process the $b$ subsets mentioned above are
counted including\vspace*{2pt} their combinations within $a+b$ sets,
which explains the binomial coefficient
$ {a+b \choose a}$ on the right-hand side.
\end{pf}

\section{Random Charlier polynomials}
\label{s5}
In order to simplify the presentation of our results
it will sometimes be convenient to use the symbolic notation
%
%
\begin{equation}
\label{defdlta2}
\Delta_{s_0} \cdots \Delta_{s_j} \prod_{p=0}^n u_{s_p}
= \mathop{\sum_{\Theta_0 \cup\cdots\cup\Theta_n = \{ 0 , 1 , \ldots, j \}
}}_{0 \notin\Theta_0 , \ldots, j \notin\Theta_j} D_{\Theta_0}
u_{s_0} \cdots D_{\Theta_n} u_{s_n} ,
\end{equation}
$s_0,\ldots,s_n \in X$, $0 \leq j \leq n$,
for any measurable process $u \dvtx \Omega^X \times X \to\real$.

The above formula implies in particular
$\Delta_{s_0} u_{s_0} = 0$, and it can be used
to rewrite the Skorohod isometry (\ref{skois}) as
\[
E_\sigma[ \delta( u )^2 ] = E_\sigma\bigl[ \Vert u \Vert^2_{L^2_\sigma(X)} \bigr]
+ E_\sigma\biggl[ \int_X \int_X \Delta_s \Delta_t ( u_t u_s ) \sigma( ds )
\sigma( dt ) \biggr] ,
\]
since by definition we have
\[
\Delta_s \Delta_t ( u_s u_t ) = D_s u_t D_t u_s,\qquad s, t \in X.
\]
In this section we show the following proposition.
\begin{prop}
\label{p12}
Let $n \geq1$ and let
$A_1(\omega) , \ldots, A_n (\omega) $ be
a.e. disjoint random Borel sets, all of them being
a.s. contained in a fixed compact set $K$ of $X$.
Then we have
\begin{eqnarray*}
&&
E_\sigma\Biggl[ \prod_{i=1}^n C_{k_i} \bigl( \delta( {\mathbf1}_{A_i} ) + \sigma( A_i
) , \sigma( A_i ) \bigr) \Biggr]
\\
&&\qquad = E_\sigma\biggl[ \int_{K^N} \Delta_{s_1} \cdots \Delta_{s_N} (
{\mathbf1}_{A_1^{k_1}} \otimes \cdots \otimes {\mathbf1}_{A_n^{k_n}} ) ( s_1,
\ldots, s_N ) \sigma( ds_1 ) \cdots\sigma( ds_N ) \biggr] ,
\end{eqnarray*}
$k_1,\ldots,k_n \in\inte$, with $N = k_1 + \cdots+ k_n$. %
\end{prop}

For $n=1$, Proposition \ref{p12} yields, in particular,
\begin{eqnarray*}
&&E_\sigma[ C_n ( \omega( A ) , \sigma( A ) ) ] \\[2pt]
&&\qquad= E_\sigma\Biggl[ \int_{K^n}
\Delta_{s_1} \cdots\Delta_{s_n} \prod_{p=1}^n {\mathbf1}_A ( s_p ) \sigma(
ds_1 ) \cdots\sigma( ds_n ) \Biggr]
\end{eqnarray*}
for $A$ a.s. contained in a fixed compact set $K$ of $X$,
which leads to (\ref{adfgsdg}) by Lem\-ma~\ref{l-4}
under condition (\ref{fdsjklfds}),
as in the following corollary which is used for the proof of
Proposition \ref{pc12.0}.
\begin{corollary}
\label{pc12}
Assume that $\tau\dvtx \Omega^X \times X \to X$ satisfies
%
%
\begin{equation}
\label{123}
D_t \tau( \omega, t ) = 0,\qquad \omega\in\Omega^X, t \in X.
\end{equation}
Then for all deterministic disjoint
$B_1,\ldots,B_n \in\mathcal{B}(X)$ we have
\begin{eqnarray*}
&&E_\sigma\Biggl[ \prod_{i=1}^n C_{k_i} \bigl( \omega\bigl( A \cap\tau^{-1} ( B_i ) \bigr) ,
\sigma\bigl( A \cap\tau^{-1} ( B_i )  \bigr) \bigr)\Biggr]
\\[2pt]
&&\qquad = E_\sigma\biggl[ \int_{A^N} D_{s_1} \cdots D_{s_N} \bigl( (
{\mathbf1}_{B_1^{k_1}} \otimes \cdots \otimes {\mathbf1}_{B_n^{k_n}} ) ( \tau(
\omega, s_1 ) , \ldots, \tau( \omega, s_N ) ) \bigr)\\[2pt]
&&\qquad\hspace*{210.1pt}{}\times \sigma( ds_1 )
\cdots\sigma( ds_N ) \biggr] ,
\end{eqnarray*}
$k_1,\ldots,k_n \in\inte$, with $N = k_1 + \cdots+ k_n$,
for all compact $A\in\mathcal{B}(X)$.
\end{corollary}
\begin{pf}
We apply Proposition \ref{p12} by letting
$A_i ( \omega) = A \cap\tau^{-1} ( \omega, B_i )$,
and we note that we have
\[
\sigma( A_i ( \omega) ) = \int_A {\mathbf1}_{B_i} ( \tau( \omega, t ) )
\sigma( dt )
= \sigma\bigl( A \cap\tau^{-1} ( \omega, B_i ) \bigr) .
\]
On the other hand, by (\ref{123}) we have
$D_t {\mathbf1}_{A_i} ( t) =
D_t {\mathbf1}_{B_i} ( \tau( \omega, t ) ) = 0
$;
hence from Lemma \ref{l0} below we have
\[
\delta( {\mathbf1}_{A_i} ) + \sigma( A_i )  = \delta( {\mathbf1}_A
{\mathbf1}_{B_i} \circ\tau) + \sigma\bigl( A \cap\tau^{-1} ( B_i ) \bigr)
=  \omega\bigl( A \cap\tau^{-1} ( B_i ) \bigr) .
\]
Finally we note that from (\ref{defdlta2}) and (\ref{123}) we have
\[
D_{s_1} \cdots D_{s_N} = \Delta_{s_1} \cdots\Delta_{s_N} ,
\]
and we apply Proposition \ref{p12}.
\end{pf}

The proof of Proposition \ref{p12} relies on the following
lemma.
\begin{lemma}
\label{l-4}
Let $F\dvtx\Omega^X \to\real$ be a bounded random variable, and
consider a random set $A$,\vadjust{\goodbreak}
a.s. contained in a fixed compact set $K$ of $X$.
For all $k \geq1$ we have
\begin{eqnarray*}
&&E_\sigma\bigl[ F C_k \bigl( \delta( {\mathbf1}_A ) + \sigma( A ) , \sigma( A ) \bigr) \bigr]
\\
&&\qquad = \sum_{z=0}^k (-1)^{k-z} \pmatrix{k \cr z} E_\sigma\Biggl[ \int_{X^k}
\prod_{j=1}^z ( I + D_{s_j} ) F \prod_{p=1}^k \mathop{\prod_{j=1}}_{j\not=
p}^z ( I + D_{s_j} ) {\mathbf1}_A ( s_p ) \\
&&\qquad\quad\hspace*{215.2pt}{}\times\sigma( ds_1 ) \cdots\sigma( ds_k
) \Biggr] .\vadjust{\goodbreak}
\end{eqnarray*}
\end{lemma}
\begin{pf}
Using Proposition \ref{l221} and Lemma \ref{ll} we have
\begin{eqnarray*}
&&E_\sigma[ F B_n ( \delta( {\mathbf1}_A ) , \sigma( A ) ) ]
\\
&&\qquad = \sum_{i=0}^n \pmatrix{n \cr i} E_\sigma[ F ( \delta( {\mathbf1}_A )
)^i B_{n-i} ( 0 , \sigma( A ) ) ]
\\
&&\qquad = \sum_{i=0}^n \pmatrix{n \cr i} \sum_{c=0}^{n-i} (-1)^c S_2 (n-i ,
c ) E_\sigma[ F ( \delta( {\mathbf1}_A ) )^i \sigma( A )^c ]
\\
&&\qquad = \sum_{i=0}^n \pmatrix{n \cr i} \sum_{c=0}^{n-i} (-1)^c S_2 (n-i ,
c) \\
&&\qquad\quad\hphantom{\sum_{i=0}^n \pmatrix{n \cr i} \sum_{c=0}^{n-i}}
\hspace*{-3pt}{}\times\sum_{e=0}^i \sum_{z=0}^e (-1)^{e-z} \pmatrix{i \cr z} S ( i - z ,
e-z )
\\
&&\qquad\quad\hspace*{92.48pt}{}\times E_\sigma\Biggl[ \int_{X^z} \Biggl( \prod_{j=1}^z ( I + D_{s_j} ) ( F \sigma( A
)^{c+e-z} ) \Biggr)\\
&&\qquad\quad\hspace*{47.6pt}\hspace*{92.2pt}{}\times \prod_{p=1}^z \mathop{\prod_{j=1}}_{j\not= p}^z ( I + D_{s_j}
) {\mathbf1}_A ( s_p ) \\
&&\qquad\quad\hspace*{184pt}{}\times\sigma( ds_1 ) \cdots\sigma( ds_z ) \Biggr]
\\
&&\qquad = \sum_{k=0}^{n-1} \sum_{i=0}^n \pmatrix{n \cr i} \sum_{c=0}^{n-i}
S_2 (n-i , c)\\
&&\qquad\quad\hspace*{69.3pt}{}\times \sum_{z=0}^{k-c} (-1)^{k-z} \pmatrix{i \cr z} S ( i - z ,
k-c-z )
\\
&&\qquad\quad\hspace*{94.8pt}{}\times E_\sigma\Biggl[ \int_{X^z} \Biggl( \prod_{j=1}^z ( I + D_{s_j} ) ( F \sigma( A
)^{k-z} ) \Biggr) \\
&&\qquad\quad\hspace*{94.8pt}\hspace*{47.6pt}{}\times\prod_{p=1}^z \mathop{\prod_{j=1}}_{j\not= p}^z ( I + D_{s_j} )
{\mathbf1}_A ( s_p ) \\[-2pt]
&&\qquad\quad\hspace*{187pt}{}\times\sigma( ds_1 ) \cdots\sigma( ds_z ) \Biggr]
\\[-2pt]
&&\qquad = \sum_{k=0}^{n-1} \sum_{z=0}^k (-1)^{k-z} \sum_{i=0}^n \pmatrix{n
\cr i} \sum_{c=0}^{n-i} \pmatrix{i \cr z} S_2 (n-i , c ) S ( i - z,
k-c-z )
\\[-2pt]
&&\qquad\quad\hspace*{29.2pt}{}\times E_\sigma\Biggl[ \int_{X^z} \Biggl( \prod_{j=1}^z ( I + D_{s_j} ) ( F
\sigma( A )^{k-z} ) \Biggr) \\[-2pt]
&&\qquad\quad\hspace*{76.7pt}{}\times\prod_{p=1}^z \mathop{\prod_{j=1}}_{j\not= p}^z ( I +
D_{s_j} ) {\mathbf1}_A ( s_p )\sigma( ds_1 ) \cdots\sigma( ds_z ) \Biggr]
\\[-2pt]
&&\qquad = \sum_{k=0}^n S ( n , k ) \sum_{z=0}^k (-1)^{k-z} \pmatrix{k \cr
z}
\\[-2pt]
&&\qquad\quad\hspace*{61.7pt}{}\times E_\sigma\Biggl[ \int_{X^z} \Biggl( \prod_{j=1}^z ( I + D_{s_j} ) ( F \sigma( A
)^{k-z} ) \Biggr) \\[-2pt]
&&\qquad\quad\hspace*{47.6pt}\hspace*{61.7pt}{}\times\prod_{p=1}^z \mathop{\prod_{j=1}}_{j\not= p}^z ( I + D_{s_j} )
{\mathbf1}_A ( s_p ) \sigma( ds_1 ) \cdots\sigma( ds_z ) \Biggr] .
\end{eqnarray*}
Hence from Proposition \ref{p1} or the
inversion formula (\ref{inv}) we get
\begin{eqnarray*}
&&
E_\sigma\bigl[ F C_k \bigl( \delta( {\mathbf1}_A ) + \sigma( A ) , \sigma( A ) \bigr) \bigr]
\\[-2pt]
&&\qquad = \sum_{z=0}^k (-1)^{k-z} \pmatrix{k \cr z}\\[-2pt]
&&\qquad\quad\hspace*{12.4pt}{}\times E_\sigma\Biggl[ \int_{X^k}
\prod_{j=1}^z ( I + D_{s_j} ) F \prod_{p=1}^k \mathop{\prod_{j=1}}_{j\not=
p}^z ( I + D_{s_j} ) {\mathbf1}_A ( s_p ) \\[-2pt]
&&\qquad\quad\hspace*{163.7pt}{}\times \sigma( ds_1 ) \cdots\sigma( ds_k
) \Biggr] .
\end{eqnarray*}
\upqed\end{pf}

In particular, Lemma \ref{l-4} applied to $F={\mathbf1}$ shows that
\begin{eqnarray*}
&&E_\sigma\bigl[ C_k \bigl( \delta( {\mathbf1}_A ) + \sigma( A ) , \sigma( A ) \bigr) \bigr]
\\
&&\qquad = \sum_{z=0}^k (-1)^{k-z} \pmatrix{k \cr z} E_\sigma\Biggl[ \int_{X^k}
\prod_{p=1}^k \mathop{\prod_{j=1}}_{j\not= p}^z ( I + D_{s_j} ) {\mathbf1}_A (
s_p ) \sigma( ds_1 ) \cdots\sigma( ds_k ) \Biggr]
\\
&&\qquad = \sum_{z=0}^k (-1)^{k-z} \pmatrix{k \cr z} E_\sigma\Biggl[ \int_{X^k}
\prod_{p=1}^k \prod_{j=1}^z ( I + \Delta_{s_j} ) {\mathbf1}_A ( s_p )
\sigma( ds_1 ) \cdots\sigma( ds_k ) \Biggr]
\\[1pt]
&&\qquad = E_\sigma\Biggl[ \int_{X^k} \Biggl( \prod_{j=1}^k ( I + \Delta_{s_j} - I ) \Biggr)
\prod_{p=1}^k {\mathbf1}_A ( s_p ) \sigma( ds_1 ) \cdots\sigma( ds_k ) \Biggr]
\\[1pt]
&&\qquad = E_\sigma\Biggl[ \int_{X^k} \Delta_{s_1} \cdots\Delta_{s_k}
\prod_{p=1}^k {\mathbf1}_A ( s_p ) \sigma( ds_1 ) \cdots\sigma( ds_k ) \Biggr] ,
\end{eqnarray*}
which is Proposition \ref{p12} for $n=1$.
Next we will apply this argument to prove
Proposition \ref{p12} from Lemma \ref{l-4}
by induction.
\begin{pf*}{Proof of Proposition \ref{p12}}
From Lemma \ref{l-4} we have
%
%
\begin{eqnarray}
\label{frl}
&&
E_\sigma\bigl[ F C_{k_1} \bigl( \delta( {\mathbf1}_{A_1} ) + \sigma( A_1
) , \sigma( A_1 ) \bigr) \bigr]
\nonumber\\[1pt]
&&\qquad = \sum_{z_1=0}^{k_1} (-1)^{k_1-z_1} \pmatrix{k_1 \cr
z_1}\nonumber\\[-8pt]\\[-8pt]
&&\qquad\quad\hspace*{16.2pt}{}\times E_\sigma
\Biggl[ \int_{K^{k_1}} \prod_{j=1}^{z_1} ( I + D_{s_{1,j}} ) F
\prod_{p=1}^{k_1} \mathop{\prod_{j=1}}_{j\not= p}^{z_1} ( I + D_{s_{1,j}} )
{\mathbf1}_{A_1} ( s_{1,p} ) \nonumber\\[1pt]
&&\qquad\quad\hspace*{178pt}{}\times\sigma( ds_{1,1} ) \cdots\sigma( ds_{1,k_1} )
\Biggr] .\nonumber
\end{eqnarray}
The first induction step is to apply the above equality to
the random variable
\[
F = \prod_{i=2}^n C_{k_i} \bigl( \delta( {\mathbf1}_{A_i} ) + \sigma( A_i ) ,
\sigma( A_i ) \bigr) .
\]
Here $F$ is not bounded, however since $A_i (\omega) \subset K$, a.s.,
$i=1,\ldots,n$, for a~fixed compact $K\in\mathcal{B}(X)$,
we check that $|F|$ is bounded by a polynomial in~$\omega(K)$, and
\[
\Biggl| \prod_{j=1}^{z_1} ( I + D_{s_{1,j}} ) F \Biggr|
\]
is bounded by another a polynomial in $\omega(K)$,
uniformly in $s_1,\ldots,s_{k_1} \in X$.
Hence by dominated convergence
we can extend (\ref{frl}) from the bounded random variable
$\max( \min( F , -C ) , C )$, $C>0$, to $F$ by letting $C$ go
to infinity.
From relation (\ref{commrel2}) we have
\begin{eqnarray*}
\prod_{j=1}^{z_1} ( I + D_{s_{1,j}} ) \delta( {\mathbf1}_{A_i} )
&=&
\delta \Biggl( \prod_{j=1}^{z_1} ( I + D_{s_{1,j}} ) {\mathbf1}_{A_i} \Biggr)
+
\sum_{k=1}^{z_1} \mathop{\prod_{j=1}}_{j \not= k}^{z_1} ( I + D_{s_{1,j}} )
{\mathbf1}_{A_i} ( s_{1,k} )
\\[-4pt]
&=& \delta \Biggl( \prod_{j=1}^{z_1} ( I + D_{s_{1,j}} ) {\mathbf1}_{A_i} \Biggr) ,
\end{eqnarray*}
$0 \leq z_1 \leq k_1$, $i \geq2$, when
$
s_{1,k} \in
{\prod^{z_1}_{j=1}}_{j \not= k}
( I + D_{s_{1,j}} )
A_1$,
$1 \leq k \leq k_1$,
since
\[
\mathop{\prod_{j=1}}_{j \not= k}^{z_1} ( I + D_{s_{1,j}} ) A_1 (\omega)
,\ldots , \mathop{\prod_{j=1}}_{j \not= k}^{z_1} ( I + D_{s_{1,j}} ) A_n
(\omega)
\]
are disjoint, $1 \leq k \leq k_1$, $\omega\in\Omega^X$, hence
\begin{eqnarray*}
\prod_{j=1}^{z_1} ( I + D_{s_{1,j}} ) F
&=& \prod_{j=1}^{z_1} ( I +
D_{s_{1,j}} )
\prod_{i=2}^n C_{k_i} \bigl( \delta( {\mathbf1}_{A_i} ) + \sigma(
A_i ) , \sigma( A_i ) \bigr)
\\[-4pt]
&=& \prod_{i=2}^n C_{k_i} \Biggl( \prod_{j=1}^{z_1} ( I + D_{s_{1,j}} )
\delta( {\mathbf1}_{A_i} ) + \prod_{j=1}^{z_1} ( I + D_{s_{1,j}} ) \sigma(
A_i ) ,\\[-4pt]
&&\hspace*{147.60pt} \prod_{j=1}^{z_1} ( I + D_{s_{1,j}} ) \sigma( A_i ) \Biggr)
\\[-4pt]
&=& \prod_{i=2}^n C_{k_i} \Biggl( \delta \Biggl( \prod_{j=1}^{z_1} ( I +
D_{s_{1,j}} ) {\mathbf1}_{A_i} \Biggr) + \prod_{j=1}^{z_1} ( I + D_{s_{1,j}} )
\sigma( A_i ) ,\\[-4pt]
&&\hspace*{152.1pt}  \prod_{j=1}^{z_1} ( I + D_{s_{1,j}} ) \sigma( A_i ) \Biggr) ,
\end{eqnarray*}
which yields, from (\ref{frl}),
\begin{eqnarray*}
\hspace*{-4pt}&&E_\sigma\Biggl[ \prod_{i=1}^n C_{k_i} \bigl( \delta( {\mathbf1}_{A_i}) + \sigma( A_i )
, \sigma( A_i ) \bigr) \Biggr]
\\[-4pt]
\hspace*{-4pt}&&\qquad = \sum_{z_1=0}^{k_1} (-1)^{k_1-z_1} \pmatrix{k_1 \cr z_1}
\\[-4pt]
\hspace*{-4pt}&&\hspace*{-4pt}\hphantom{\sum_{z_1=0}^{k_1}}\qquad\quad{}
\times E_\sigma \Biggl[ \int_{X^{k_1}} \prod_{i=2}^n C_{k_i} \Biggl( \delta \Biggl(
\prod_{j=1}^{z_1} ( I + D_{s_{1,j}} ) {\mathbf1}_{A_i} \Biggr)
+\prod_{j=1}^{z_1} ( I + D_{s_{1,j}} ) \sigma( A_i ) ,\\[-4pt]
\hspace*{-4pt}&&\hspace*{-6pt}\qquad\quad\hspace*{226pt} \prod_{j=1}^{z_1}
( I + D_{s_{1,j}} ) \sigma( A_i ) \Biggr)
\\[-4pt]
\hspace*{-4pt}&&\qquad\quad\hspace*{86pt}{}\times \prod_{p=1}^{k_1} \mathop{\prod_{j=1}}_{j\not= p}^{z_1} ( I +
D_{s_{1,j}} ) {\mathbf1}_{A_1} ( s_{1,p} ) \sigma( ds_{1,1} ) \cdots\sigma(
ds_{1,k_1} ) \Biggr] .
\end{eqnarray*}
Next, we apply Lemma \ref{l-4} again to
\begin{eqnarray*}
&&C_{k_2} \Biggl( \delta \Biggl( \prod_{j=1}^{z_1} ( I + D_{s_{1,j}} ) {\mathbf1}_{A_2} \Biggr)
\\[-3pt]
&&\qquad\hspace*{0pt}{}+ \prod_{j=1}^{z_1} ( I + D_{s_{1,j}} ) \sigma( A_2 ) ,
\prod_{j=1}^{z_1} ( I + D_{s_{1,j}} ) \sigma( A_2 ) \Biggr)
\end{eqnarray*}
and to
\begin{eqnarray*}
F &=& \prod_{p=1}^{k_1} \mathop{\prod_{j=1}}_{j\not= p}^{z_1} ( I +
D_{s_{1,j}} ) {\mathbf1}_{A_1} ( s_{1,p} )
\\[-3pt]
& &\hphantom{\prod_{p=1}^{k_1} \mathop{\prod_{j=1}}_{j\not= p}^{z_1} }
\hspace*{0pt}{} \times \prod_{i=3}^n C_{k_i} \Biggl( \delta \Biggl( \prod_{j=1}^{z_1} ( I +
D_{s_{1,j}} ) {\mathbf1}_{A_i} \Biggr) + \prod_{j=1}^{z_1} ( I + D_{s_{1,j}} )
\sigma( A_i ) ,\\[-3pt]
&&\qquad\quad\hspace*{166.6pt} \prod_{j=1}^{z_1} ( I + D_{s_{1,j}} ) \sigma( A_i ) \Biggr) ,
\end{eqnarray*}
and by iteration of this argument we obtain
\begin{eqnarray*}
\hspace*{-6pt}&&E_\sigma\Biggl[ \prod_{i=1}^n C_{k_i} \bigl( \delta( {\mathbf1}_{A_i}) + \sigma( A_i )
, \sigma( A_i ) \bigr) \Biggr]
\\[-3pt]
\hspace*{-6pt}&&\qquad = \sum_{z_n=0}^{k_n} \cdots \sum_{z_1=0}^{k_1} \prod_{l=1}^n
(-1)^{k_l-z_l} \\[-3pt]
\hspace*{-6pt}&&\qquad\quad\hspace*{49.7pt}{}\times \prod_{l=1}^n \pmatrix{k_l \cr z_l }
E_\sigma \Biggl[ \int_{X^{k_1}} \prod_{i=1}^n \mathop{\prod_{j=1}}_{j\not=
p}^{z_i} ( I + D_{s_{i,j}} ) \\[-3pt]
\hspace*{-6pt}&&\qquad\quad\hspace*{141.pt}{}\times \prod_{i=1}^n \prod_{j=1}^{k_i}
{\mathbf1}_{A_i} ( s_{i,j} ) \sigma( ds_{1,1} ) \cdots\sigma( ds_{n,k_n} )
\Biggr]
\\[-3pt]
\hspace*{-6pt}&&\qquad = \sum_{z_n=0}^{k_n} \cdots \sum_{z_1=0}^{k_1} \prod_{l=1}^n
(-1)^{k_l-z_l} \\[-3pt]
\hspace*{-6pt}&&\qquad\quad\hspace*{49.7pt}{}\times \prod_{l=1}^n \pmatrix{k_l \cr z_l }E_\sigma
\Biggl[ \int_{X^{k_1}} \Biggl( \prod_{i=1}^n \prod_{j=1}^{z_i} ( I +
\Delta_{s_{i,j}} ) \Biggr) \\[-3pt]
\hspace*{-6pt}&&\qquad\quad\hspace*{142.1pt}{}\times\prod_{i=1}^n \prod_{j=1}^{k_i} {\mathbf1}_{A_i} (
s_{i,j} ) \sigma( ds_{1,1} ) \cdots\sigma( ds_{n,k_n} ) \Biggr]
\\[-3pt]
\hspace*{-6pt}&&\qquad = E_\sigma\Biggl[ \int_{X^N} \Biggl( \prod_{i=1}^n \prod_{j=1}^{k_i} ( I +
\Delta_{s_{i,j}} - I ) \Biggr)\\[-3pt]
\hspace*{-6pt}&&\qquad\quad\hspace*{37.6pt}{}\times \prod_{i=1}^n \prod_{j=1}^{k_i} {\mathbf1}_{A_i} (
s_{i,j} ) \sigma( ds_{1,1} ) \cdots\sigma( ds_{n,k_n} ) \Biggr]
\\[-3pt]
\hspace*{-6pt}&&\qquad = E_\sigma\Biggl[ \int_{X^N} \Biggl( \prod_{i=1}^n \prod_{j=1}^{k_i}
\Delta_{s_{i,j}} \Biggr) \prod_{i=1}^n \prod_{j=1}^{k_i} {\mathbf1}_{A_i} (
s_{i,j} ) \sigma( ds_{1,1} ) \cdots\sigma( ds_{n,k_n} ) \Biggr] .
\end{eqnarray*}
\upqed\end{pf*}

Next we prove Proposition \ref{t1.1}.
\begin{pf*}{Proof of Proposition \ref{t1.1}}
Taking $g \dvtx Y \to\real$ to be the step function
\[
g = \sum_{i=1}^m c_i {\mathbf1}_{B_i} ,
\]
where $c_1,\ldots,c_m \in\real$ and $B_1,\ldots, B_m \in\mathcal
{B} ( Y )$
are disjoint Borel subsets of $Y$, Corollary \ref{pc12} shows that for
compact $A \in\mathcal{B}(X)$ we have
\begin{eqnarray*}
&&
E_\sigma\biggl[ e^{ - \int_A g( \tau( \omega, t ) ) \sigma(dt) } \prod_{x \in
A \cap\omega} \bigl( 1 + g ( \tau( \omega, x ) ) \bigr) \biggr] \\
&&\qquad= E_\sigma\Biggl[
\prod_{l=1}^m e^{ - c_l \sigma( A \cap\tau^{-1} ( B_l ) ) }
\prod_{l=1}^m ( 1 + c_l )^{ \omega( A \cap\tau^{-1} ( B_l ) )} \Biggr]
\\
&&\qquad = \sum_{k_1=0}^\infty \cdots \sum_{k_m=0}^\infty \Biggl( \prod_{i=1}^m
\frac{c_i^{k_i}}{k_i!} \Biggr) E_\sigma\Biggl[ \prod_{i=1}^m C_{k_i} \bigl( \omega\bigl( A
\cap\tau^{-1} ( B_i ) \bigr) , \sigma\bigl( A \cap\tau ^{-1} ( B_i ) \bigr) \bigr)\Biggr]
\\
&&\qquad = \sum_{n=0}^\infty \frac{1}{n!} E_\sigma[ I_n ( {\mathbf1}_{A^n} (
\cdot) g^{\otimes n} ( \tau^{\otimes n} ( \omega, \cdot) ) ) ]
\\
&&\qquad = \sum_{n=0}^\infty \frac{1}{n!} E_\sigma\Biggl[ \int_{A^n} \Delta_{s_1}
\cdots\Delta_{s_n} \prod_{p=1}^n g ( \tau( \omega, s_p ) ) \sigma( ds_1
) \cdots\sigma( ds_n ) \Biggr] .
\end{eqnarray*}
In the general case with $g \dvtx Y \to\real$ bounded measurable
the conclusion follows by approximation of $g$
by step functions and dominated convergence under (\ref{st1}),
followed by extension to $A=X$ using the bound (\ref{bnd1.1}).
\end{pf*}

Finally we state the following lemma which has been
used in the proof of Corollary \ref{pc12}.
\begin{lemma}
\label{l0}
Assume that
%
%
\begin{equation}
\label{xy}
D_t \tau( \omega, t ) = 0,\qquad \omega\in\Omega^X, t \in X.
\end{equation}
Then we have
\[
\int_X {\mathbf1}_A ( t ) h ( \tau( \omega, t ) ) \omega( d t ) = \delta(
{\mathbf1}_A h \circ\tau) ( \omega) + \int_A h \circ\tau( \omega, t )
\sigma( dt ),\qquad \omega\in\Omega^X,
\]
for all compact $A \in\mathcal{B}(X)$ and
all bounded measurable functions $h\dvtx X \to\real$.
\end{lemma}
\begin{pf}
We note that condition (\ref{xy}) above means that
$\tau(\omega,t)$ does not depend on the presence
or absence of a point in $\omega$ at $t$, and in particular,
\[
\tau( \omega, t ) = \tau(\omega\cup\{ t\} , t ) ,\qquad t \notin\omega ,
\]
and
\[
\tau( \omega, t ) = \tau(\omega\setminus\{ t\} , t ) ,\qquad t \in\omega .
\]
Hence we have
\begin{eqnarray*}
&&\delta( {\mathbf1}_A h \circ\tau) + \int_A h \circ\tau( \omega, t ) \sigma(
dt )
\\
&&\qquad = \int_X {\mathbf1}_A ( t ) h \bigl( \tau( \omega\setminus\{ t \} , t ) \bigr) \bigl(
\omega( d t ) - \sigma( dt ) \bigr) + \int_X {\mathbf1}_A ( t ) h ( \tau(
\omega, t ) ) \sigma( dt )
\\
&&\qquad = \int_X {\mathbf1}_A ( t ) h \bigl( \tau( \omega\setminus\{ t \} , t ) \bigr)
\omega( d t )
\\
&&\qquad = \int_X {\mathbf1}_A ( t ) h ( \tau( \omega, t ) ) \omega( d t ) .
\end{eqnarray*}
\upqed\end{pf}

\section{Link with the Carleman--Fredholm determinant}
\label{s6}
In this section we make some remarks on
differences between the Poisson and Wiener cases,
in relation to the quasi-nilpotence of random transformations.
We consider a Poisson random measure on
$\real_+\times[-1,1]^d$ on the real line with flat intensity measure,
in which case it is known \cite{girpri,girunif,priqi},
that,
building the Poisson measure as a product of exponential and
uniform densities on the sequence space $\real^\inte$, we have the
Girsanov identity,
\[
E [ F ( I+u) | {\det}_2 (I+\nabla u) | \exp(- \nabla^* (u)) ] = E [ F ]
,
\]
where $u\dvtx \real^\inte\to\real^\inte$ is a random
shift satisfying certain conditions,
${\det}_2 (I+\nabla u)$
is the Carleman--Fredholm determinant of
$I+\nabla u$ and $\nabla^* (u)$ is a~Skoro\-hod-type
integral of the discrete-time process $u$.

When it is invertible,
$(I+u)_*\pi_\sigma$
is absolutely continuous with respect to~$\pi_\sigma$ with
\[
\frac{d(I+u)^{-1}_* \pi_\sigma}{d\pi_\sigma} = | {\det}_2 (I+\nabla u)
| \exp(-\nabla^* (u)) .
\]
It can be checked (cf. \cite{girpri,girunif,priqi})
that in the adapted case this
yields the usual Girsanov theorem for the change of intensity
of Poisson random measures when the configuration points are shifted
by an adapted smooth diffeomorphism $\phi\dvtx \Omega^X \times
\real_+\times[0,1]^d
\longrightarrow\real_+\times[0,1]^d$,
in which case $I+Du$ becomes a block diagonal matrix, each $d\times d$
block having the Jacobian determinant
$\vert\partial_{t,x} \phi(\omega, T_k,x_k^1,\ldots, x_k^d) \vert$,
and we have
\begin{eqnarray*}
&&{\det}_2 (I +\nabla u)\exp(-\nabla^* ( u )) \\
&&\qquad= e^{
-\int_{\real_+\times[0,1]^{d}} (\vert\partial_{s,x} \phi( \omega, s,x)
\vert-1 ) \,ds\,dx } \prod_{k=1}^\infty\vert\partial_{t,x} \phi(\omega,
T_k,x_k^1,\ldots, x_k^d ) \vert.
\end{eqnarray*}
The main difference with the Wiener case is that
here $\nabla u$ is not quasi-nilpo\-tent on $\ell^2 ( \inte)$ and
we do not have ${\det}_2 (I +\nabla u) = 1$.
Nevertheless it should be possible to recover Proposition \ref{pc12.0}
in a weaker form by checking the relation
\[
{\det} (I +\nabla u) = \prod_{k=1}^\infty\vert\partial_{t,x}
\phi(\omega, T_k,x_k^1,\ldots, x_k^d ) \vert
\]
for anticipating shifts $\phi\dvtx \Omega^X \times\real_+\times[0,1]^d
\longrightarrow\real_+\times[0,1]^d$,
under smoothness and quasi-nilpotence assumptions stronger than those
assumed in this paper.


%

%
\printaddresses

\end{document}